\documentclass[a4paper,12pt, reqno]{amsart}
\usepackage{txfonts}      
\usepackage{amssymb}
\usepackage{eucal}
\usepackage{amsmath}
\usepackage{amscd}
\usepackage[dvips]{color}
\usepackage{mathtools}
\usepackage{multicol}
\usepackage[all]{xy}           
\usepackage{graphicx}
\usepackage{color}
\usepackage{colordvi}
\usepackage{xspace}
\usepackage[normalem]{ulem}
\usepackage{cancel}
\usepackage{tikz}
\usepackage{longtable}
\usepackage{multirow}
\usepackage{ifpdf}
\ifpdf
\usepackage[colorlinks,final,backref=page,hyperindex]{hyperref}
\else
\usepackage[colorlinks,final,backref=page,hyperindex,hypertex]{hyperref}
\fi
\usepackage{placeins}

\begin{document}
\newcommand {\emptycomment}[1]{} 

\newcommand{\tabincell}[2]{\begin{tabular}{@{}#1@{}}#2\end{tabular}}

\newcommand{\nc}{\newcommand}
\newcommand{\delete}[1]{}


\delete{
\nc{\label}[1]{\label{#1}  
{\hfill \hspace{1cm}{\bf{{\ }\hfill(#1)}}}}
\nc{\cite}[1]{\cite{#1}{{\bf{{\ }(#1)}}}}  
\nc{\ref}[1]{\ref{#1}{{\bf{{\ }(#1)}}}}  
\nc{\eqref}[1]{~\eqref{#1}{{\bf{{\ }(#1)}}}} 
\nc{\bibitem}[1]{\bibitem[\bf #1]{#1}} 
}

\newtheorem{thm}{Theorem}[section]
\newtheorem{lem}[thm]{Lemma}
\newtheorem{cor}[thm]{Corollary}
\newtheorem{pro}[thm]{Proposition}
\newtheorem{conj}[thm]{Conjecture}
\theoremstyle{definition}
\newtheorem{defi}[thm]{Definition}
\newtheorem{ex}[thm]{Example}
\newtheorem{rmk}[thm]{Remark}
\newtheorem{pdef}[thm]{Proposition-Definition}
\newtheorem{condition}[thm]{Condition}

\renewcommand{\labelenumi}{{\rm(\alph{enumi})}}
\renewcommand{\theenumi}{\alph{enumi}}
\renewcommand{\labelenumii}{{\rm(\roman{enumii})}}
\renewcommand{\theenumii}{\roman{enumii}}

\nc{\tred}[1]{\textcolor{red}{#1}}
\nc{\tblue}[1]{\textcolor{blue}{#1}}
\nc{\tgreen}[1]{\textcolor{green}{#1}}
\nc{\tpurple}[1]{\textcolor{purple}{#1}}
\nc{\btred}[1]{\textcolor{red}{\bf #1}}
\nc{\btblue}[1]{\textcolor{blue}{\bf #1}}
\nc{\btgreen}[1]{\textcolor{green}{\bf #1}}
\nc{\btpurple}[1]{\textcolor{purple}{\bf #1}}


\newcommand{\End}{\text{End}}

\nc{\calb}{\mathcal{B}}
\nc{\call}{{L}}
\nc{\calo}{\mathcal{O}}
\nc{\frakg}{A}
\nc{\frakh}{\mathfrak{h}}
\nc{\ad}{\mathrm{ad}}

\nc{\ccred}[1]{\tred{\textcircled{#1}}}


\newcommand{\cm}[1]{\textcolor{red}{\underline{CM:}#1 }}
\newcommand{\gl}[1]{\textcolor{blue}{\underline{GL:}#1 }}
\nc{\li}[1]{\textcolor{red}{#1}}
\nc{\lir}[1]{\textcolor{red}{\underline{Li:}#1}}

\title[A unified approach to Manin triples of Lie algebras]
{A unified approach to Manin triples of Lie algebras and the corresponding bialgebra structures}

\author{Guilai Liu}
\address{Chern Institute of Mathematics \& LPMC, Nankai University, Tianjin 300071, China
    \newline
    Institut de Recherche Math\'ematique Avanc\'ee,
    UMR 7501, Universit\'e de Strasbourg et CNRS,
    7 rue Ren\'e-Descartes, 67000 Strasbourg, France
}
\email{liugl@mail.nankai.edu.cn}

\author{Chengming Bai}
\address{Chern Institute of Mathematics \& LPMC, Nankai University, Tianjin 300071, China}
\email{baicm@nankai.edu.cn}

\author{Li Guo}
\address{Department of Mathematics and Computer Science, Rutgers University, Newark, NJ 07102, USA}
\email{liguo@rutgers.edu}

\begin{abstract}
    We introduce self-module Lie bialgebras as a common framework for bialgebra structures associated with different types of Manin triples of Lie algebras. A self-module Lie algebra is a Lie algebra with a
    specified representation on its underlying vector space; equivalently, it is an NS-Lie algebra expressed in
    different notation. We define a generic Manin triple as a finite-dimensional Lie algebra equipped with a nondegenerate bilinear form and two complementary isotropic Lie subalgebras, and establish the equivalence among standard generic Manin triples, suitable matched pairs of Lie algebras, families of matched pairs of self-module Lie algebras and self-module Lie bialgebras. The
    framework recovers Lie, pre-Lie and anti-pre-Lie bialgebras as special cases and admits examples
    outside these three families.
    We also introduce the self-module Lie Yang–Baxter equation and show
    that its solutions yield self-module Lie bialgebras and admit the $\mathcal{O}$-operator characterization. Finally,
    we obtain solutions in semi-direct product self-module Lie algebras, including constructions arising
    from self-bimodule pre-Lie algebras.
\end{abstract}

\subjclass[2020]{
17B62,  
17B38,  
16T10,	
17A36,  
17B10, 
17B60, 
17D25.  
}

\keywords{representation of Lie algebra, self-module Lie algebra, Manin triple, bialgebra, Yang-Baxter equation}

\maketitle

\vspace{-1cm}

\tableofcontents

\vspace{-1cm}

\allowdisplaybreaks

\section{Introduction}
Lie bialgebras \cite{Cha,Drinfeld1983} play a central role in the theories of quantum groups, Yang-Baxter equations, and are closely related to Poisson-Lie groups.
A remarkable property of Lie bialgebras is their description by Manin triples of Lie algebras equipped with nondegenerate
symmetric invariant bilinear forms. 

This characterization has been generalized to other algebraic structures, including pre-Lie bialgebras and anti-pre-Lie bialgebras. The pre-Lie bialgebra \cite{Bai2008} was introduced in the study of parak\"ahler Lie algebras and admits an equivalent characterization as a Manin triple of Lie algebras with respect to the symplectic form. The anti-pre-Lie bialgebra was introduced in \cite{LB2023} as a fundamental ingredient in the bialgebra theory for transposed Poisson algebras \cite{TPA}, which is equivalently characterized as a Manin triple of Lie algebras with respect to the nondegenerate commutative $2$-cocycle.

Despite their different origins, their Manin triple characterizations have the following properties.
 \begin{enumerate}
    \item\label{r1-} There is a Lie algebra $(\frak g\oplus\frak h,[-,-]_{d})$ containing
    $(\frak g,[-,-]_{\frak g})$ and $(\frak h,[-,-]_{\frak h})$ as Lie subalgebras, and a nondegenerate bilinear form $\mathcal{B}_{d} $ on $\frak g\oplus \frak h$ such that
    $ \frak g $ and $ \frak h $ are isotropic with respect to $\mathcal{B}_{d} $.
    \item\label{r2-}  The bilinear form $\mathcal{B}_{d}$ is either symmetric or antisymmetric.
    \item\label{r3-} There are additional compatibility conditions of $\mathcal{B}_{d}$ with the Lie bracket $[-,-]_d$ on $\mathfrak g\oplus\mathfrak h$.
 \end{enumerate}

These observations suggest a unified framework. We introduce the notion of a {\bf generic Manin triple} of Lie algebras by retaining only property \eqref{r1-} that is shared by the Manin triples of Lie algebras corresponding to Lie bialgebras, pre-Lie bialgebras and anti-pre-Lie bialgebras, while dropping the requirements \eqref{r2-} and \eqref{r3-} in which these Manin triples exhibit different forms.

To study the structure of a generic Manin triple of Lie algebras, we first examine the underlying structure of a Lie algebra with a nondegenerate bilinear form.
For this purpose, in Section~\ref{sec:2}, we introduce the notion of a self-module Lie algebra
$(\frak g,[-,-]_{\frak g},\circ_{\frak g})$,
 consisting of a Lie algebra $(\frak g,[-,-]_{\frak g})$
 together with a multiplication
 $\circ_{\frak g}:\frak g\otimes\frak g\rightarrow\frak g$
 subject to an additional condition ensuring that
 $(L_{\circ_{\frak g}},\frak g)$ is a representation of
 $(\frak g,[-,-]_{\frak g})$ (Definition~\ref{def:2.1}).
This notion is an equivalent reformulation of NS-Lie
 algebras introduced in \cite{SR-Lie}. We use the name {\bf self-module} to emphasize that the Lie algebra $\mathfrak{g}$ acts on its own underlying vector space via $x\mapsto L_{\circ_{\mathfrak{g}}}(x)$. This general notion will be applied to 
 Lie algebras, pre-Lie algebras \cite{Bai2021.2,Bur} and
 anti-pre-Lie algebras \cite{LB2022} to recover their established bialgebra theories, as well as to develop new bialgebra structures. 

A Lie algebra $(\frak g,[-,-]_{\frak g})$ with a nondegenerate
bilinear form $\mathcal{B}$ induces a self-module Lie algebra
structure $(\frak g,[-,-]_{\frak g},\circ_{\frak g})$. 
An additional linear isomorphism $\phi$ on $\frak g$ also yields a self-module Lie algebra $(\frak g,[-,-]_{\frak g},\circ_{\phi})$, called the $\phi$-conjugate self-module Lie algebra of $(\frak g,[-,-]_{\frak g},\circ_{\frak g})$.
To generalize the matched pairs of Lie algebras as a device connecting Manin triples of Lie algebras and Lie bialgebras, we study representations and matched pairs of self-module Lie algebras (Definition~\ref{defi:mpSR-Lie}). 

In Section \ref{sec:3}, we introduce generic Manin triples of Lie algebras (Definition~\ref{991}) and use them to extract the conditions for self-module Lie bialgebras.
We first show that every finite-dimensional generic Manin triple of Lie algebras is isomorphic to a standard one on $\frak g\oplus\frak g^*$ equipped with the nondegenerate bilinear form $\mathcal B_{\phi}$ for some linear isomorphism $\phi$ on $\frak g$ (Proposition \ref{1005}), where
\begin{eqnarray}\label{eq:intro-Bphi1}
    \mathcal{B}_{\phi}(x+\xi ,y+\eta )=\langle \phi(x),\eta \rangle+\langle \xi ,y\rangle,\;\forall x,y\in\frak g, \xi ,\eta \in\frak g^{*}.
\end{eqnarray}
The underlying Lie bracket is determined by a matched pair of Lie algebras 
$\big((\frak g,[-,-]_{\frak g}),(\frak g^*$,
$[-,-]_{\frak g^*}),L^*_{\circ_{\frak g}},L^*_{\circ_{\frak g^*}}\big)$,
in which the multiplications $\circ_{\frak g}$ and $\circ_{\frak g^*}$ make $(\frak g,[-,-]_{\frak g},\circ_{\frak g})$ and $(\frak g^*,[-,-]_{\frak g^*}$,
$\circ_{\frak g^*})$ into self-module Lie algebras respectively (Lemma \ref{1452}).
The self-module Lie algebra $(\frak g\oplus\frak g^*,[-,-]_d,\bullet_{\phi})$
induced from $\mathcal B_{\phi}$ gives rise to a family of matched pairs of
self-module Lie algebras (Proposition \ref{1503}) by varying $\phi$.

Dualizing the Lie bracket and self-module product on $\frak g^{*}$ yields two comultiplications on $\frak g$. Together with the self-module Lie algebra structure and suitable compatibility conditions, these define a self-module Lie bialgebra.
As one of the main results in this paper, we establish an equivalence among a standard generic Manin triple of Lie algebras, a matched pair of Lie algebras, a family of matched pairs of self-module Lie algebras and a self-module Lie bialgebra (Theorem \ref{cor:2532}).
The notion of self-module Lie bialgebras is sufficiently broad to encompass all the bialgebra structures admitting an equivalent description in terms of generic Manin triples of Lie algebras. In particular, self-module Lie bialgebras provide a unified framework for the study of such bialgebra structures, including the known Lie, pre-Lie and anti-pre-Lie bialgebras.

In Section \ref{sec:4}, we introduce the notion of a self-module Lie Yang-Baxter equation (SLYBE) in a self-module Lie algebra, which plays a role analogous to that of the classical Yang-Baxter equation in the theory of Lie bialgebras.
 Every solution $r\in\frak g\otimes\frak g$ of the SLYBE gives rise to a self-module Lie bialgebra (Corollary \ref{cor:4.3}).
 If $r$ moreover satisfies the $\phi$-condition for a linear isomorphism $\phi$ on $\frak g$, then the resulting self-module Lie bialgebra is triangular with respect to $\phi$.
 Triangular Lie bialgebras and triangular anti-pre-Lie bialgebras are recovered by taking $\phi=\mathrm{id}_{\frak g}$, whereas triangular pre-Lie bialgebras are recovered by taking $\phi=-\mathrm{id}_{\frak g}$.
Solutions of the SLYBE are characterized by $\mathcal{O}$-operators of the underlying Lie algebra and,
under the $\phi$-condition, of a $\phi$-conjugate self-module Lie algebra (Theorem \ref{thm:3034}).
Self-bimodule pre-Lie algebras are used to construct $\mathcal{O}$-operators of self-module Lie algebras and hence solutions of the SLYBE in the corresponding semi-direct product self-module Lie algebras (Proposition \ref{pro:pre-Ns-} and Corollary \ref{2110}).
Example \ref{ex:nonclassical} provides a triangular self-module Lie bialgebra that shows that our general framework goes beyond the three classical cases. 

The main correspondences and constructions for self-module Lie bialgebras are illustrated in Figure~\ref{fig:main-correspondences} on page \pageref{fig:main-correspondences}.
Finally, Section \ref{2584} records the three classical bialgebra specializations. 

\smallskip

\noindent
{\bf Notation.} Throughout this paper, unless otherwise specified, all vector spaces are assumed to be
 finite-dimensional over a field $\mathbb{K}$ of characteristic $0$, over which all the linear maps and tensor products are taken.
 For a vector space $\frak g$ with a multiplication
 $\circ_{\frak g}:\frak g\otimes\frak g\to\frak g$,
 we define $L_{\circ_{\frak g}}(x)y=x\circ_{\frak g}y=R_{\circ_{\frak g}}(y)x$ for all $x,y\in\frak g$.
In particular, if $(\frak g,[-,-]_{\frak g})$ is a Lie algebra, we write
 $\mathrm{ad}_{\frak g}(x)y=[x,y]_{\frak g}$.
 We denote the identity map by $\mathrm{id}_{\frak g}$ or $1$ if there is no risk of confusion,
 and use $\tau$ for the flip map on $\frak g\otimes\frak g$ defined by
 $\tau(x\otimes y)=y\otimes x$ for all $x,y\in \frak g$.

\vspace{-.35cm}

\section{Self-module Lie algebras and their matched pairs}\label{sec:2}

We study properties and matched pairs of self-module Lie algebras. We show that
a Lie algebra together with a nondegenerate bilinear form gives rise to a self-module Lie algebra.
A self-module Lie algebra with a linear isomorphism of the underlying space gives rise to a conjugate self-module Lie algebra. We also study representations and matched pairs of self-module Lie algebras in preparation for the study of generic Manin triples of Lie algebras in Section~\ref{sec:3}.

\subsection{Self-module Lie algebras and nondegenerate bilinear forms on Lie algebras}\

\begin{defi}\label{def:2.1}
A {\bf self-module Lie algebra} is a triple $(\frak g,[-,-]_{\frak g},\circ_{\frak g})$ such that $(\frak g,[-,-]_{\frak g})$ is a Lie algebra, $\circ_{\frak g}:\frak g\otimes\frak g\rightarrow\frak g$ is a multiplication  and the following equation holds:
\begin{equation}\label{468}
[x,y]_{\frak g}\circ_{\frak g}z=x\circ_{\frak g}(y\circ_{\frak g}z)-y\circ_{\frak g}(x\circ_{\frak g}z),\;\forall x,y,z\in\frak g.
\end{equation}
\end{defi}

Recall that a {\bf representation} of a Lie algebra $(\frak
g,[-,-]_{\frak g})$ is a pair $(\rho,V)$ such that $V$ is a vector
space, $\rho:\frak g\rightarrow\mathrm{End}_{\mathbb K}(V)$ is a
linear map and the following equation holds:
\begin{equation*}\label{eq:repLie}
\rho([x,y]_{\frak g})v=\rho(x)\rho(y)v-\rho(y)\rho(x)v,\;\forall x,y\in \frak g, v\in V.
\end{equation*}
In fact, the pair $(\rho,V)$ is a representation of $(\frak g,[-,-]_{\frak g})$ if and only if there is a Lie algebra structure on the direct sum $\frak g\oplus V$ of vector spaces
given by
\begin{eqnarray}\label{eq:Lie sd}
[x+u,y+v]_{d}=[x,y]_{\frak g}+\rho(x)v-\rho(y)u,\;\forall x,y\in\frak g,u,v\in V.
\end{eqnarray}
In this case,
we call $(\frak g\oplus V,[-,-]_{d} )$ the {\bf semi-direct product Lie algebra of $(\frak g,[-,-]_{\frak g})$ with respect to $(\rho,V)$} and simply denote it by $ \frak g \ltimes_{\rho }V$. The following characterization justifies the terminology
of a self-module Lie algebra.

\begin{pro}\label{473}
Let $(\frak g,[-,-]_{\frak g})$ be a Lie algebra. Then there exists a multiplication $\circ_{\frak g}:\frak g\otimes \frak g\rightarrow \frak g$ such that $(\frak g,[-,-]_{\frak g},\circ_{\frak g})$ is a self-module Lie algebra if and only if there is  a representation $(\rho,\frak g)$ of  $(\frak g,[-,-]_{\frak g})$.
Here the mutual relation is given by $\rho=L_{\circ_{\frak g}}$, or equivalently,
\begin{equation*}
x\circ_{\frak g}y=\rho(x)y,\;\forall x,y\in\frak g.
\end{equation*}
\end{pro}
\begin{proof}
It follows by
a straightforward verification.
\end{proof}

\begin{rmk}\label{rmk530}
Recall \cite{SR-Lie} that an {NS-Lie algebra}  is a triple $(\frak g,\circ_{\frak g},\curlyvee_{\frak g})$ consisting of a vector space $\frak g$  and multiplications $\circ_{\frak g},\curlyvee_{\frak g}:\frak g\otimes\frak g\rightarrow\frak g$ such that \eqref{468} and  the following equations hold:
\begin{align*}
&        x\curlyvee_{\frak g}y+y\curlyvee_{\frak g}x=0,\\
&        x\curlyvee_{\frak g}[y,z]_{\frak g}+y\curlyvee_{\frak g}[z,x]_{\frak g}+z\curlyvee_{\frak g}[x,y]_{\frak g}
+x\circ_{\frak g}(y\curlyvee_{\frak g}z)
+y\circ_{\frak g}(z\curlyvee_{\frak g}x)
+z\circ_{\frak g}(x\curlyvee_{\frak g}y)=0,
\end{align*}
for all $x,y,z\in\frak g$, where $[-,-]_{\frak g}:\frak g\otimes\frak g\rightarrow\frak g$ is given by
\begin{eqnarray}\label{eq:SR-Lie3}
[x,y]_{\frak g}=x\circ_{\frak g}y-y\circ_{\frak g}x+x\curlyvee_{\frak g}y,\;\forall x,y\in\frak g.
\end{eqnarray}
Suppose that $(\frak g,\circ_{\frak g},\curlyvee_{\frak g})$ is an NS-Lie algebra. Then $(\frak g,[-,-]_{\frak g},\circ_{\frak g})$ with $[-,-]_{\frak g} $ given by \eqref{eq:SR-Lie3} is a self-module Lie algebra.
Conversely, every self-module Lie algebra
$(\mathfrak{g},[-,-]_{\mathfrak{g}},\circ_{\mathfrak{g}})$
determines an NS-Lie algebra
$(\mathfrak{g},\circ_{\mathfrak{g}},\curlyvee_{\mathfrak{g}})$
with $\curlyvee_{\mathfrak{g}}$ defined through \eqref{eq:SR-Lie3}.
Consequently, there is a one-to-one correspondence between self-module Lie algebras and NS-Lie algebras.
\end{rmk}

Lie algebras, pre-Lie algebras and anti-pre-Lie algebras naturally
provide self-module Lie algebra structures (Example \ref{3234} and Propositions \ref{2515}, \ref{2780}).

Recall that two representations $(\mu,U)$ and  $( \rho ,V )$ of a Lie algebra $(\frak g,[-,-]_{\frak g})$
are called {\bf equivalent} if there exists a linear isomorphism $f:U\rightarrow V$ such that
\begin{eqnarray}\label{eq:434}
f\big(\mu(x)u\big )=\rho(x)f(u),\;\forall x\in\frak g,  u \in U.
\end{eqnarray}

Let $\frak g,V$ be vector spaces and $\rho:\frak g\rightarrow\mathrm{End}_{\mathbb K}(V)$ be a linear map.
Define a linear map $\rho^{*}:\frak g\rightarrow \mathrm{End}_{\mathbb K}(V^{*})$ by
\begin{eqnarray*}
    \langle \rho^{*}(x)u^{*},v\rangle=-\langle u^{*},\rho(x)v\rangle,\;\forall x\in \frak g, u^{*}\in V^{*},v\in V,
\end{eqnarray*}
where $\langle-,-\rangle$ denotes the natural duality pairings.

\begin{pro}\label{pro:comp}
Let $(\frak g,[-,-]_{\frak g} )$ be a Lie algebra and $\mathcal{B}$ be a nondegenerate bilinear form on $\frak g$.  Then there is a self-module Lie algebra $(\frak g,[-,-]_{\frak g},\circ_{\frak g})$ in which the multiplication $\circ_{\frak g}:\frak g\otimes\frak g\rightarrow \frak g$ is given by
\begin{equation}
\mathcal{B}\big(x\circ_{\frak g}y,z\big)=-\mathcal{B}(y,[x,z]_{\frak g}),\;\forall x,y,z\in \frak g.\label{eq:pro:comp1}
\end{equation}
In this case, we say $(\frak g,[-,-]_{\frak g},\circ_{\frak g})$ is {\bf the self-module Lie algebra induced from $(\frak g,[-,-]_{\frak g},\mathcal{B})$}.
\end{pro}
\begin{proof}
Note that \eqref{eq:pro:comp1} equivalently translates into the following equation:
\begin{equation*}\label{eq:823}
L_{\circ_{\frak g}}(x)=(\mathcal{B}^{\natural})^{-1} \mathrm{ad}_{\frak g}^{*}(x)\mathcal{B}^{\natural},\;\forall x\in\frak g,
\end{equation*}
where  $\mathcal{B}^{\natural}:\frak g\rightarrow \frak g^{*}$ is the linear isomorphism given by
\begin{equation*}
    \langle\mathcal{B}^{\natural}(x),y\rangle=\mathcal{B}(x,y),\;\forall x,y\in\frak g.
\end{equation*}
Thus
$(L_{\circ_{\frak g}},\frak g)$ is a representation of $(\frak g,[-,-]_{\frak g})$ which is equivalent to the coadjoint representation $(\mathrm{ad}^{*}_{\frak g},\frak g^{*})$. Therefore $(\frak g,[-,-]_{\frak g},\circ_{\frak g})$ is a self-module Lie algebra.
\end{proof}

Let $\frak g,U,V$ be vector spaces and
$\phi:\frak g\rightarrow\frak g$,
$f:V\rightarrow U$, $g:U\rightarrow V$,
$\rho:\frak g\rightarrow\mathrm{End}_{\mathbb K}(V)$ be linear maps.
Define linear maps
$\rho[f,\phi,g]:\frak g\rightarrow\mathrm{End}_{\mathbb K}(U)$
and $\phi+f:\frak g\oplus V\rightarrow\frak g\oplus U$ by
\begin{equation*}
\rho[f,\phi,g](x)u
 =f\rho\big(\phi(x)\big)g(u),\;
 (\phi+f)(x+v)
 =\phi(x)+f(v),
 \quad\forall x\in\frak g,u\in U,v\in V.
\end{equation*}

\begin{lem}\label{lem:2.7}
Given a self-module Lie algebra $(\frak g,[-,-]_{\frak g},\circ_{\frak g})$ and a linear isomorphism $\phi:\frak g\rightarrow \frak g$,
there is a self-module Lie algebra
$(\frak g,[-,-]_{\frak g},\circ_{\phi})$ with product  $\circ_{\phi}$ given by
\begin{equation}\label{612}
x\circ_{\phi}y=\phi^{-1}\big( x\circ_{\frak g}\phi(y)  \big),\;\forall x,y\in\frak g.
\end{equation}
In this case, we say $(\frak g,[-,-]_{\frak g},\circ_{\phi})$
is the {\bf $\phi$-conjugate self-module Lie algebra  of $(\frak g,[-,-]_{\frak g},\circ_{\frak g})$},
or $(\frak g,[-,-]_{\frak g},\circ_{\phi})$ is {\bf conjugate to $(\frak g,[-,-]_{\frak g},\circ_{\frak g})$ through $\phi$}.
Moreover, the resulting semi-direct product Lie algebras $ \frak g\ltimes_{L_{\circ_{\frak g}}}\frak g$ and $ \frak g\ltimes_{L_{\circ_{\phi}}}\frak g$ are isomorphic through $\mathrm{id}_{\frak g}+\phi^{-1}$.
\end{lem}
\begin{proof}
By Proposition \ref{473}, $(L_{\circ_{\frak g}},\frak g)$ is a representation of $(\frak g,[-,-]_{\frak g})$.
Since $\phi:\frak g\rightarrow \frak g$ is a linear isomorphism,
 the pair $ (  {L}_{\circ_{\frak g}}[\phi^{-1},1,\phi],\frak g\big)$ is a representation of $(\frak g,[-,-]_{\frak g})$ that is equivalent to  $(L_{\circ_{\frak g}},\frak g)$. By Proposition \ref{473} again, there is a self-module Lie algebra $(\frak g,[-,-]_{\frak g},\circ_{\phi})$ with the 
multiplication  $\circ_{\phi}$  given by
\begin{equation*}
x\circ_{\phi}y={L}_{\circ_{\frak g}}[\phi^{-1},1,\phi](x)y=\phi^{-1} L_{\circ_{\frak g}} (x)\phi(y)=\phi^{-1}\big( x\circ_{\frak g}\phi(y)  \big),\;\forall x,y\in\frak g.
\end{equation*}
The remaining result can be checked directly.
\end{proof}

\begin{pro}\label{1493-}
Let $(\frak g,[-,-]_{\frak g})$ be a Lie algebra, $\mathcal{B}$ be a nondegenerate bilinear form on $\frak g$ and $(\frak g,[-,-]_{\frak g}$,
$\circ_{\frak g})$ be the  self-module Lie algebra induced from $(\frak g,[-,-]_{\frak g},\mathcal{B})$ through Proposition \ref{pro:comp}.
Suppose that $(\frak g,[-,-]_{\frak g})$ admits another self-module Lie algebra structure $(\frak g,[-,-]_{\frak g},\star_{\frak g})$.
Then there exists a linear isomorphism $\phi$ such that $(\frak g,[-,-]_{\frak g},\star_{\frak g})$ is conjugate to $(\frak g,[-,-]_{\frak g},\circ_{\frak g})$ through $\phi$ if and only if there exists a nondegenerate bilinear form $\omega$ on $\frak g$ such that $(\frak g,[-,-]_{\frak g},\star_{\frak g})$ is the self-module Lie algebra induced from $(\frak g,[-,-]_{\frak g},\omega)$.
Here the mutual relation is given by
\begin{equation}\label{1498}
\omega(x,y)=\mathcal{B}\big(\phi(x),y\big),\;\forall x,y\in\frak g.
\end{equation}
\end{pro}
\begin{proof}
Suppose that there exists a linear isomorphism $\phi$ such that $(\frak g,[-,-]_{\frak g},\star_{\frak g})$ is conjugate to $(\frak g,[-,-]_{\frak g}$,
$\circ_{\frak g})$ through $\phi$.
For all $x,y,z\in\frak g$, we have
\begin{equation*}
\omega(x\star_{\frak g}y,z)=\omega\Big(  \phi^{-1}\big(x\circ_{\frak g}\phi(y)\big),z\Big)\overset{\eqref{1498}}{=}\mathcal{B}\big(  x\circ_{\frak g}\phi(y),z \big)=-\mathcal{B}\big( \phi(y),[x,z]_{\frak g} \big)\overset{\eqref{1498}}{=}-\omega(y,[x,z]_{\frak g}).
\end{equation*}
Hence $(\frak g,[-,-]_{\frak g},\star_{\frak g})$ is the self-module Lie algebra induced from $(\frak g,[-,-]_{\frak g},\omega)$.
Conversely, suppose that there exists a nondegenerate bilinear form $\omega$ on $\frak g$ such that $(\frak g,[-,-]_{\frak g},\star_{\frak g})$ is the self-module Lie algebra induced from $(\frak g,[-,-]_{\frak g},\omega)$. In a similar way, we can show that $(\frak g,[-,-]_{\frak g},\star_{\frak g})$ is conjugate to $(\frak g,[-,-]_{\frak g},\circ_{\frak g})$ through the linear isomorphism $\phi:\frak g\rightarrow\frak g$ given through \eqref{1498}.
\end{proof}

\begin{ex}\label{557}
Let $(\frak n=\mathrm{span}_{\mathbb{K}}\lbrace e_{1}, e_{2}\rbrace,[-,-]_{\frak n})$ be the $2$-dimensional nonabelian Lie algebra whose nonzero product is
\begin{equation}\label{572}
 [e_{1}, e_{2}]_{\frak n}=e_{1}.
\end{equation}
Let $\omega$ be the nondegenerate bilinear form on $\frak n$ with the nonzero values given by
\begin{equation*}
\omega(e_{1},e_{2})=1=-\omega(e_{2}, e_{1}).
\end{equation*}
Then there is a self-module Lie algebra $(\frak n,[-,-]_{\frak n},\circ_{\frak n})$ given in Proposition \ref{pro:comp}, where the nonzero products of $\circ_{\frak n}$ are given by
\begin{equation}\label{566}
e_{1}\circ_{\frak n}e_{2}=e_{1},\;
e_{2}\circ_{\frak n}e_{2}=e_{2}.
\end{equation}
Moreover, by Proposition \ref{1493-}, any self-module Lie algebra structure on $(\frak n,[-,-]_{\frak n})$ induced from a nondegenerate bilinear form on $\frak n$ is conjugate to $(\frak n,[-,-]_{\frak n},\circ_{\frak n})$ through some linear isomorphism on $\frak n$.
\end{ex}

\subsection{Representations and matched pairs of self-module Lie algebras}

\begin{defi}
Let $(\frak g,[-,-]_{\frak g},\circ_{\frak g})$ be a self-module
Lie algebra and $(\rho,V),(l,V)$ be representations of $(\frak
g,[-,-]_{\frak g})$. If there is a linear map $r:\frak
g\rightarrow\mathrm{End}_{\mathbb K}(V)$ such that the following
equation holds:
\begin{equation}\label{eq:867}
r(y)\rho(x)v-l(x)r(y)v +r (x\circ_{\frak g}y )v=0,\;\forall x,y\in\frak g,v\in V,
\end{equation}
then $(\rho,l,r,V)$ is called a {\bf representation} of $(\frak g,[-,-]_{\frak g},\circ_{\frak g})$.
Two representations $(\mu,s,t,U)$ and $(\rho,l,r,V)$ of a self-module Lie algebra $(\frak g,[-,-]_{\frak g},\circ_{\frak g})$ are called {\bf equivalent}
if there exists a linear isomorphism $f:U\rightarrow V$ such that \eqref{eq:434} and the following equations hold:
\begin{equation*}\label{eq:684}
f\big(s(x)u\big )=l(x)f(u),\;
f\big(t(x)u\big )=r(x)f(u),
\;\forall x\in\frak g, u \in U.
\end{equation*}
\end{defi}

\begin{pro}\label{pro:2.14}
Let $(\frak g,[-,-]_{\frak g},\circ_{\frak g})$ be a self-module
Lie algebra, $V$ be a vector space and $\rho,l,r:\frak
g\rightarrow\mathrm{End}_{\mathbb K}(V)$ be linear maps. Then
$(\rho,l,r, V)$ is a representation of $(\frak g,[-,-]_{\frak
g},\circ_{\frak g})$ if and only if there is  a self-module Lie
algebra structure $ ( \frak g\oplus V,[-,-]_{d},\circ_{d} )$ on
$\frak g\oplus V$ with $[-,-]_d$ given by \eqref{eq:Lie sd} and
$\circ_d$ given by
\begin{eqnarray}
(x+u)\circ_{d}  (y+v)=x\circ_{\frak g}y+ l(x)v+r(y)u, \;\forall x,y\in \frak g, u,v\in V.\label{693}
\end{eqnarray}
We call  $(\frak g\oplus V,[-,-]_{d},\circ_{d})$ the {\bf semi-direct product self-module Lie algebra
of $(\frak g,[-,-]_{\frak g},\circ_{\frak g})$ with respect to $(\rho,l,r, V)$} and simply denote it by $\frak g
\ltimes_{\rho,l,r}V$.
\end{pro}
\begin{proof}
    The result follows from Proposition \ref{pro:mat} by taking $\frak h=V$ with the zero Lie bracket and multiplication and setting $\rho_{\frak h}=l_{\frak h}=r_{\frak h}=0$.
\end{proof}

\begin{ex}
Let $(\frak g,[-,-]_{\frak g},\circ_{\frak g})$ be a self-module Lie algebra. Then $(\mathrm{ad}_{\frak g},L_{\circ_{\frak g}}, R_{\circ_{\frak g}},\frak g)$ is a representation of
$(\frak g,[-,-]_{\frak g},\circ_{\frak g})$,  called the {\bf adjoint representation} of $(\frak g,[-,-]_{\frak g},\circ_{\frak g})$.
\end{ex}

\begin{pro}\label{752}
Let $\frak g=(\frak g,[-,-]_{\frak g},\circ_{\frak g})$ be a self-module Lie algebra,
$\phi:\frak g\rightarrow\frak g$ be a linear isomorphism and $\frak g_{\phi}=(\frak g,[-,-]_{\frak g},\circ_{\phi})$ be the $\phi$-conjugate self-module Lie algebra of $(\frak g,[-,-]_{\frak g},\circ_{\frak g})$.
Suppose that  $(\rho,l,r,V)$ is a representation of $(\frak g,[-,-]_{\frak g},\circ_{\frak g})$ and $f:V\rightarrow V$ is a linear isomorphism.
There is a self-module Lie algebra $\frak g_{\phi}\ltimes_{ \rho, l[f^{-1},1, f], r[f^{-1},\phi,1 ] }V$ which is conjugate to $\frak g\ltimes_{\rho,l,r}V$ through $\phi+f$.
Consequently,
the quadruple
$ ( \rho,l[f^{-1},1,  f], r[f^{-1},\phi,1] ,V )$   is a representation of $(\frak g,[-,-]_{\frak g},\circ_{\phi})$.
\end{pro}
\begin{proof}
By the assumption, it follows from Proposition \ref{pro:2.14} that there is a self-module Lie algebra structure $(\frak g\oplus V,[-,-]_{d},\circ_{d})=\frak g\ltimes_{\rho,l,r}V$ given by \eqref{eq:Lie sd} and \eqref{693}.
Since $\phi+f$ is a linear isomorphism on $\frak g\oplus V$, there is a self-module Lie algebra $(\frak g\oplus V,[-,-]_{d},\circ'_{d})$ which is conjugate to  $(\frak g\oplus V,[-,-]_{d},\circ_{d})$ through $\phi+f$, in which the multiplication $\circ'_{d}$ is explicitly given by
\begin{align*}
(x+u)\circ'_{d}(y+v)&= (\phi+f )^{-1}\Big( (x+u)\circ_{d}(\phi+f)(y+v)  \Big)\\
&\overset{\eqref{693}}{=} (\phi^{-1}+f^{-1})\Big( x\circ_{\frak g}\phi(y)+l(x)f(v)+r\big(
\phi(y)
\big)u \Big)\\
&= \phi^{-1}\big( x\circ_{\frak g}\phi(y) \big)+f^{-1}l(x)f(v)+f^{-1}  r\big(\phi (y) \big) u  \\
&\overset{\eqref{612}}{=} x\circ_{\phi}y+l[f^{-1},1, f](x)v+r[f^{-1},\phi ,1](y)u,\;\forall x,y\in\frak g, u,v\in V.
\end{align*}
By Proposition \ref{pro:2.14}, $(\rho, l[f^{-1},1,  f], r[f^{-1},\phi,1] ,V)$
is a representation of $(\frak g,[-,-]_{\frak g},\circ_{\phi})$.
\end{proof}

\begin{lem}\label{801n}
	Let $(\frak g,[-,-]_{\frak g},\circ_{\frak g})$ be a self-module Lie algebra and $(\rho,l,r,V)$ be a representation of $(\frak g,[-,-]_{\frak g},\circ_{\frak g})$.
	Then $(l^{*},\rho^{*},r^{*},V^{*})$ is a representation of $(\frak g,[-,-]_{\frak g},\circ_{\frak g})$.
\end{lem}
\begin{proof}
	It is clear that $(l^{*},V^{*})$ and $(\rho^{*},V^{*})$ are representations of $(\frak g,[-,-]_{\frak g})$. Moreover, for all $x,y\in\frak g, u^{*}\in V^{*},v\in V$, we have
	\begin{equation*}
		\langle r^{*}(y)l^{*}(x)u^{*}-\rho^{*}(x)r^{*}(y)u^{*}+r^{*}(x\circ_{\frak g}y)u^{*},v\rangle=\langle u^{*}, l(x)r(y)v-r(y)\rho(x)v-r(x\circ_{\frak g}y)v\rangle\overset{\eqref{eq:867}}{=}0,
	\end{equation*}
	that is,
	\begin{equation*}
		r^{*}(y)l^{*}(x)u^{*}-\rho^{*}(x)r^{*}(y)u^{*}+r^{*}(x\circ_{\frak g}y)u^{*}=0.
	\end{equation*}
	Therefore, $(l^{*},\rho^{*},r^{*},V^{*})$ is a representation of $(\frak g,[-,-]_{\frak g},\circ_{\frak g})$.
\end{proof}

Let $U$ and $V$ be vector spaces and
$f:U\rightarrow V$ be a linear map.
Define its transpose $f^{t}:V^{*}\rightarrow U^{*}$ by
\begin{equation*}
    \langle f^{t}(\eta), u\rangle
    =\langle \eta,f(u )\rangle,
    \;\forall u\in U, \eta\in V^{*}.
\end{equation*}

\begin{pro}\label{848}
Under the assumptions of Proposition \ref{752}, the quadruple $$\big( l^{*},\rho^{*}[(f^{-1})^{t},1,f^{t}],r^{*}[(f^{-1})^{t},\phi,1],V^{*} \big)$$ is a representation of $(\frak g,[-,-]_{\frak g},\circ_{\phi})$. Moreover, the resulting semi-direct product self-module Lie algebra $\frak g_{\phi}\ltimes_{ l^{*},\rho^{*}[(f^{-1})^{t},1,f^{t}],r^{*}[(f^{-1})^{t},\phi,1] }V^{*}$ is conjugate to $\frak g\ltimes_{l^{*},\rho^{*},r^{*}}V^{*}$ through $\phi+f^{t}$.
\end{pro}
\begin{proof}
By Proposition \ref{752}, $(\rho, l[f^{-1},1,  f], r[f^{-1},\phi,1] ,V)$
is a representation of the self-module Lie algebra $(\frak g,[-,-]_{\frak g},\circ_{\phi})$.
By a straightforward verification, the quadruple
\begin{align*}
&\big(  \rho[f,1,f^{-1}], (l[f^{-1},1,  f])[f,1,f^{-1}], (r[f^{-1},\phi,1])  [f,1,f^{-1}],V \big)\\
=&(  \rho[f,1,f^{-1}],l,  r[1,\phi,f^{-1}], V )
\end{align*}
is a representation of $(\frak g,[-,-]_{\frak g},\circ_{\phi})$,
which is equivalent to
$( \rho, l[f^{-1},1,  f], r[f^{-1},\phi,1] ,V )$.
By Lemma \ref{801n}, the quadruple
\begin{equation*}
\big( l^{*}, (\rho[f,1,f^{-1}])^{*}, (r[1,\phi,f^{-1}])^{*}, V^{*} \big)
=\big( l^{*},\rho^{*}[(f^{-1})^{t},1,f^{t}],r^{*}[(f^{-1})^{t},\phi,1],V^{*} \big)
\end{equation*}
is a representation of $(\frak g,[-,-]_{\frak g},\circ_{\phi})$.
The remaining result follows from a straightforward verification.
\end{proof}

\begin{pro}\label{776}
	Let $\frak g=(\frak g,[-,-]_{\frak g},\circ_{\frak g})$ be a self-module Lie algebra, $\phi,\psi:\frak g\rightarrow\frak g$ be arbitrary linear isomorphisms and $\frak g_{\phi}=(\frak g,[-,-]_{\frak g},\circ_{\phi})$ be the  $\phi$-conjugate self-module Lie algebra of $(\frak g,[-,-]_{\frak g},\circ_{\frak g})$.
Then $(L^{*}_{\circ_{\frak g}},\mathrm{ad}^{*}_{\frak g}, R^{*}_{\circ_{\frak g}} ,{\frak g}^* )$ is a representation of $(\frak g,[-,-]_{\frak g},\circ_{\frak g})$ and   $ (L^{*}_{\circ_{\frak g }},\mathrm{ad}^{*}_{\frak g}[(\psi^{-1})^{t},1$,
$\psi^{t}],R^{*}_{\circ_{\frak g}}[(\psi^{-1})^{t},\phi,1],\frak g^{*})$  is a representation of $(\frak g,[-,-]_{\frak g},\circ_{\phi})$.
	Moreover, the resulting self-module Lie algebra  $\frak g_{\phi}\ltimes_{ L^{*}_{\circ_{\frak g }},\mathrm{ad}^{*}_{\frak g}[(\psi^{-1})^{t}, 1, \psi^{t}],
		R^{*}_{\circ_{\frak g}}[(\psi^{-1})^{t},\phi,1]  }\frak g^{*}$
	 is conjugate to $\frak g\ltimes_{{L}_{\circ_{\frak g}}^{*},\mathrm{ad}^{*}_{\frak g}, R^{*}_{\circ_{\frak g}} }\frak g^{*}$ through $\phi+\psi^{t}$.
\end{pro}
\begin{proof}
	It follows from Propositions \ref{752} and \ref{848} by taking $(\rho,l,r,V)=(\mathrm{ad}_{\frak g},L_{\circ_{\frak g}}, R_{\circ_{\frak g}},\frak g)$ and $f=\psi$. 
\end{proof}

We now recall the notion of matched pairs of Lie algebras.
\begin{defi}\cite{Maj}\label{1493}
Let $(\frak g,[-,-]_{\frak g})$ and $(\frak h,[-,-]_{\frak h})$ be
Lie algebras and $\rho_{\frak g}:\frak
g\rightarrow\mathrm{End}_{\mathbb K}(\frak h),\rho_{\frak h}:\frak
h\rightarrow\mathrm{End}_{\mathbb K}(\frak g)$ be linear maps. If
there is a Lie algebra structure on the vector space $\frak
g\oplus \frak h$ given by
\begin{equation}\label{eq:Lie-}
[x+a,y+b]_{d}=[x,y]_{\frak g}+\rho_{\frak h}(a)y-\rho_{\frak h}(b)x+[a,b]_{\frak h}+\rho_{\frak g}(x)b-\rho_{\frak g}(y)a,\;\;\forall x,y\in \frak g, a,b\in \frak h,
\end{equation}
then we say  $\big((\frak g,[-,-]_{\frak g}),(\frak h,[-,-]_{\frak
h}),\rho_{\frak g},\rho_{\frak h}\big)$ is a {\bf matched pair of
Lie algebras}. In this case, we
denote the Lie algebra structure on $\frak g\oplus\frak h$ by
$\frak g\bowtie_{\rho_{\frak g}}^{\rho_{\frak h}}\frak h$.
\end{defi}

\begin{pro}\cite{Maj}\label{1509}
Let $(\frak g,[-,-]_{\frak g})$ and $(\frak h,[-,-]_{\frak h})$ be
Lie algebras and $\rho_{\frak g}:\frak
g\rightarrow\mathrm{End}_{\mathbb K}(\frak h),\rho_{\frak h}:\frak
h\rightarrow\mathrm{End}_{\mathbb K}(\frak g)$ be linear maps.
Then $\big((\frak g,[-,-]_{\frak g}),(\frak h,[-,-]_{\frak
h}),\rho_{\frak g},\rho_{\frak h}\big)$ is a  matched pair of Lie
algebras if and only if $(\rho_{\frak g},\frak h)$ is a
representation of  $(\frak g,[-,-]_{\frak g})$, $(\rho_{\frak
h},\frak g)$ is a representation of $(\frak h,[-,-]_{\frak h})$
and the following equations hold for all $x,y\in \frak g,a,b\in
\frak h$:
\begin{align}
\rho_{\frak g}(x)[a,b]_{\frak h}-[\rho_{\frak g}(x)a,b]_{\frak h}-[a,\rho_{\frak g}(x)b]_{\frak h}+\rho_{\frak g}(\rho_{\frak h}(a)x)b-\rho_{\frak g}(\rho_{\frak h}(b)x)a&=0, \label{eq:matched pair of Lie algebras1-}\\
\rho_{\frak h}(a)[x,y]_{\frak g}-[\rho_{\frak h}(a)x,y]_{\frak g}-[x,\rho_{\frak h}(a)y]_{\frak g}+\rho_{\frak h}(\rho_{\frak g}(x)a)y-\rho_{\frak h}(\rho_{\frak g}(y)a)x&=0.\label{eq:matched pair of Lie algebras2-}
\end{align}
\end{pro}

Now we introduce the notion of matched pairs of self-module Lie algebras.
\begin{defi}\label{defi:mpSR-Lie}
Let $(\frak g,[-,-]_{\frak g},\circ_{\frak g})$ and $(\frak
h,[-,-]_{\frak h},\circ_{\frak h})$ be self-module Lie algebras
and $\rho_{\frak g},l_{\frak g},r_{\frak g}:\frak
g\rightarrow\mathrm{End}_{\mathbb K}(\frak h)$, $ \rho_{\frak
h},l_{\frak h},r_{\frak h}:\frak h\rightarrow\mathrm{End}_{\mathbb
K}(\frak g)$ be linear maps. If there is a self-module Lie algebra
$ (\frak g\oplus\frak h,[-,-]_{d}, \circ_{d} )$   with $[-,-]_{d}$
given by \eqref{eq:Lie-} and $\circ_{d}$ given by
\begin{equation*}
(x+a)\circ_{d} (y+b)=x\circ_{\frak g}y+l_{\frak h}(a)y+r_{\frak h}(b)x+a\circ_{\frak h}b+l_{\frak g}(x)b+r_{\frak g}(y)a,\;\forall x,y\in \frak g, a,b\in \frak h,
\end{equation*}
then we say the collection
\begin{equation}\label{875}
\big( (\frak g,[-,-]_{\frak g},\circ_{\frak g}),(\frak h,[-,-]_{\frak h},\circ_{\frak h} ), \rho_{\frak g},l_{\frak g}, r_{\frak g}, \rho_{\frak h},l_{\frak h}, r_{\frak h} \big)
\end{equation}
is a {\bf matched pair of self-module Lie algebras}. In this case,
we   denote the self-module Lie
algebra structure on $\frak g\oplus\frak h$ by $ \frak g
\bowtie_{\rho_{\frak g},l_{\frak g},r_{\frak g}}^{\rho_{\frak
h},l_{\frak h},r_{\frak h}} \frak h $.
\end{defi}

\begin{pro}\label{pro:mat}
Let $(\frak g,[-,-]_{\frak g},\circ_{\frak g})$ and $(\frak
h,[-,-]_{\frak h},\circ_{\frak h})$ be self-module Lie algebras
and $\rho_{\frak g},l_{\frak g},r_{\frak g}:\frak
g\rightarrow\mathrm{End}_{\mathbb K}(\frak h)$, $ \rho_{\frak
h},l_{\frak h},r_{\frak h}:\frak h\rightarrow\mathrm{End}_{\mathbb
K}(\frak g)$ be linear maps. Then the collection \eqref{875} is a
matched pair of self-module Lie algebras if and only if the
following conditions hold.
\begin{enumerate}
\item $(\rho_{\frak g},l_{\frak g},r_{\frak g},\frak h)$ is a
representation of $(\frak g,[-,-]_{\frak g},\circ_{\frak g})$.
\item $(\rho_{\frak h},l_{\frak h},r_{\frak h},\frak g)$ is a
representation of $(\frak h,[-,-]_{\frak h},\circ_{\frak h})$.
\item $\big((\frak g,[-,-]_{\frak g}),(\frak h$, $[-,-]_{\frak
h}),\rho_{\frak g},\rho_{\frak h}\big)$ is a matched pair of Lie
algebras, that is, \eqref{eq:matched pair of Lie algebras1-} and
\eqref{eq:matched pair of Lie algebras2-} hold. \item The
following equations hold for all $x,y\in \frak g, a,b\in \frak h$:
\begin{eqnarray}
r_{\frak h}(a)[x,y]_{\frak g}-x\circ_{\frak g}\big(r_{\frak h}(a)y\big)-r_{\frak h}\big( l_{\frak g}(y)a\big)x+y\circ_{\frak g}\big(r_{\frak h}(a)x\big)+r_{\frak h}\big( l_{\frak g}(x)a\big)y=0,&\label{12641}\\
l_{\frak h}\big( \rho_{\frak g}(x)a  \big)y
- \big(\rho_{\frak h}(a)x\big)\circ_{\frak g}y
-x\circ_{\frak g}\big(l_{\frak h}(a)y\big)
-r_{\frak h}\big(r_{\frak g}(y)a\big)x
+l_{\frak h}(a)(x\circ_{\frak g}y)=0,&\label{12642}\\
r_{\frak g}(x)[a,b]_{\frak h}
-a\circ_{\frak h} \big(r_{\frak g}(x)b\big)
-r_{\frak g}\big( l_{\frak h}(b)x\big)a
+b\circ_{\frak h}\big(r_{\frak g}(x)a\big)+r_{\frak g}\big( l_{\frak h}(a)x\big)b=0,&\label{12643}\\
l_{\frak g}\big( \rho_{\frak h}(a)x  \big)b
- \big( \rho_{\frak g}(x)a \big)\circ_{\frak h}b
-a\circ_{\frak h}\big(l_{\frak g}(x)b\big)
-r_{\frak g}\big(r_{\frak h}(b)x\big)a
+l_{\frak g}(x)(a\circ_{\frak h}b)=0.& \label{12644}
\end{eqnarray}
\end{enumerate}
\end{pro}
\begin{proof}
It follows from a straightforward computation.
\end{proof}

\begin{pro}\label{935}
Let $\big( (\frak g,[-,-]_{\frak g},\circ_{\frak g} ),(\frak h,[-,-]_{\frak h},\circ_{\frak h}), \rho_{\frak g},l_{\frak g}, r_{\frak g}, \rho_{\frak h},l_{\frak h}, r_{\frak h} \big)$ be a matched pair of self-module Lie algebras.
Let $\phi:\frak g\rightarrow\frak g$ and $f:\frak h\rightarrow\frak h$ be linear isomorphisms, $\frak g_{\phi}=(\frak g,[-,-]_{\frak g},\circ_{\phi})$ be the $\phi$-conjugate self-module Lie algebra of $\frak g=(\frak g,[-,-]_{\frak g},\circ_{\frak g})$ and
$\frak h_{f}=(\frak h,[-,-]_{\frak h},\circ_{f})$ be the $f$-conjugate self-module Lie algebra of $\frak h=(\frak h,[-,-]_{\frak h},\circ_{\frak h})$.
\begin{enumerate}
\item\label{939} There is a self-module Lie algebra $\frak g_{\phi}\bowtie_{\rho_{\frak g},l_{\frak g}[f^{-1},1,f],r_{\frak g}[f^{-1},\phi,1]}^{\rho_{\frak h},l_{\frak h}[\phi^{-1},1,\phi],r_{\frak h}[\phi^{-1},f,1]} \frak h_{f}$  conjugate to $\frak g \bowtie_{\rho_{\frak g},l_{\frak g},r_{\frak g}}^{\rho_{\frak h},l_{\frak h},r_{\frak h}} \frak h $ through $\phi+f$.
\item\label{940} The collection
\begin{equation*}
\big( (\frak g,[-,-]_{\frak g},\circ_{\phi}), (\frak h,[-,-]_{\frak h},\circ_{f}), \rho_{\frak g},l_{\frak g}[f^{-1},1,f],r_{\frak g}[f^{-1},\phi,1],
\rho_{\frak h},l_{\frak h}[\phi^{-1},1,\phi],r_{\frak h}[\phi^{-1},f,1]
\big)\\
\end{equation*}
is a matched pair of self-module Lie algebras.
\end{enumerate}
\end{pro}
\begin{proof}
By the assumption, it follows from Lemma \ref{lem:2.7} that there
is a self-module Lie algebra $\big(\frak g \bowtie_{\rho_{\frak g}
}^{\rho_{\frak h} } \frak h=(\frak g\oplus\frak h,[-,-]_{d}),\circ_{\phi+f}\big)$
that is $(\phi+f)$-conjugate to
$\frak g \bowtie_{\rho_{\frak g},l_{\frak g},r_{\frak
g}}^{\rho_{\frak h},l_{\frak h},r_{\frak h}} \frak h $, in which
the multiplication $\circ_{\phi+f}$ is given by  
\begin{align*}
&   (x+a)\circ_{\phi+f}(y+b)\\
&= (\phi+f)^{-1}\big( (x+a)\circ_{d}(\phi+f)(y+b)  \big)\\
&= (\phi^{-1}+f^{-1})\Big( x\circ_{\frak g}\phi(y)+l_{\frak h}(a)\phi(y)+r_{\frak h}\big( f(b)\big)x+a\circ_{\frak h}f(b)+l_{\frak g}(x)f(b)+r_{\frak g}\big(\phi(y) \big)a  \Big)\\
&= x\circ_{\phi}y+l_{\frak h}[\phi^{-1},1,\phi](a)y+ r_{\frak h}[\phi^{-1},f,1](b)x
+a\circ_{f}b+l_{\frak g}[f^{-1},1,f](x)b+  r_{\frak g}[f^{-1},\phi,1](y)a,
\end{align*}
for all $x,y\in\frak g, a,b\in\frak h$. Thus we have
\begin{equation*}
(\frak g \bowtie_{\rho_{\frak g} }^{\rho_{\frak h} } \frak h,\circ_{\phi+f})=    \frak g _{\phi} \bowtie_{\rho_{\frak g},l_{\frak g}[f^{-1},1,f],r_{\frak g}[f^{-1},\phi,1]}^{\rho_{\frak h},l_{\frak h}[\phi^{-1},1,\phi],r_{\frak h}[\phi^{-1},f,1]} \frak h _{f}.
\end{equation*}
Therefore we obtain Item \eqref{939}. Item  \eqref{940} follows from Definition \ref{defi:mpSR-Lie}.
\end{proof}

\section{Generic Manin triples of Lie algebras and self-module Lie bialgebras}\label{sec:3}

We introduce generic Manin triples of Lie algebras and show that every such triple is isomorphic to a standard one.
We then introduce self-module Lie bialgebras and establish the equivalence among standard generic Manin triples of Lie algebras, suitable matched pairs of Lie algebras, families of matched pairs of self-module Lie algebras and self-module Lie bialgebras.

\subsection{Generic Manin triples of Lie algebras}

\begin{defi}\label{991}
    Let $(\frak g,[-,-]_{\frak g})$ and
    $(\frak h,[-,-]_{\frak h})$ be Lie algebras, and let
    $(\frak g\oplus\frak h,[-,-]_d)$ be a Lie algebra containing
    them as Lie subalgebras. If there is a
     nondegenerate bilinear
        form
             $\mathcal B_d$
        on $\frak g\oplus\frak h$ such that both
        $\frak g$ and $\frak h$ are isotropic with respect to
        $\mathcal B_d$, then we call
        \[
        \big(
        (\frak g\oplus\frak h,[-,-]_d,\mathcal B_d),
        (\frak g,[-,-]_{\frak g}),
        (\frak h,[-,-]_{\frak h})
        \big)
        \]
        a {\bf generic Manin triple of Lie algebras}.
        A generic Manin triple of Lie algebras is called
        {\bf standard} if $\frak h=\frak g^*$ and
        $\mathcal B_d=\mathcal B_{\phi}$ for some linear
        isomorphism $\phi:\frak g\rightarrow\frak g$, where
        \[
        \mathcal B_{\phi}(x+\xi ,y+\eta )
        =
        \langle\phi(x),\eta \rangle+\langle \xi ,y\rangle,
        \;\forall
        x,y\in\frak g,\ \xi ,\eta \in\frak g^*.
        \]
\end{defi}

Two generic Manin triples of Lie algebras
$\big((\frak g\oplus\frak h,[-,-]_d,\mathcal B_d),(\frak g,[-,-]_{\frak g}),(\frak h,[-,-]_{\frak h})\big)$
and
$
\big(
(\frak g'\oplus\frak h',[-,-]'_d,\mathcal B'_d),
(\frak g',[-,-]'_{\frak g}),
(\frak h',[-,-]'_{\frak h})
\big)
$
are called {\bf isomorphic} if there is an isomorphism
$
\Phi:(\frak g\oplus\frak h,[-,-]_d)
\rightarrow
(\frak g'\oplus\frak h',[-,-]'_d)
$
of Lie algebras such that
\begin{equation*}
    \Phi(\frak g)=\frak g',\;\Phi(\frak h)=\frak h',\;
    \mathcal B_d(u,v)=  \mathcal B'_d(\Phi(u),\Phi(v))
    ,
    \;
    \forall u,v\in\frak g\oplus\frak h.
\end{equation*}

\begin{pro}\label{1005}
    Every generic Manin triple of Lie algebras is isomorphic to a standard one.
\end{pro}
\begin{proof}
    Suppose that
    $\big((\frak g\oplus\frak h,[-,-]_d,\mathcal B_d),
    (\frak g,[-,-]_{\frak g}),
    (\frak h,[-,-]_{\frak h})\big)$
    is a generic Manin triple of Lie algebras.
    Define a linear map $f:\frak h\rightarrow\frak g^*$ by
    $\langle f(a),x\rangle=\mathcal B_d(a,x)$ for all
    $x\in\frak g$ and $a\in\frak h$.
    Since $\mathcal B_d$ is nondegenerate and both summands
    are isotropic, $f$ is a linear isomorphism.
    Define a multiplication $[-,-]'_d$ on
    $\frak g\oplus\frak g^*$ by
    \[
    [x+\xi ,y+\eta ]'_d
    =(\mathrm{id}_{\frak g}+f)
    \big([x+f^{-1}(\xi ),y+f^{-1}(\eta )]_d\big),
    \quad
    \forall x,y\in\frak g,\ \xi ,\eta \in\frak g^*.
    \]
    Then $(\frak g\oplus\frak g^*,[-,-]'_d)$ is a Lie algebra
    containing $(\frak g,[-,-]_{\frak g})$ and
    $(\frak g^*,[-,-]_{\frak g^*}
    =[-,-]'_d|_{\frak g^*\otimes\frak g^*})$
    as Lie subalgebras.
    Similarly, there is a linear isomorphism
    $\phi:\frak g\rightarrow\frak g$ determined by
    $\langle\phi(x),f(a)\rangle=\mathcal B_d(x,a)$
    for all $x\in\frak g$ and $a\in\frak h$.
    Moreover, for all $x,y\in\frak g$ and $a,b\in\frak h$,
    we have
    \begin{equation*}
        \mathcal B_d(x+a,y+b)
 =\langle\phi(x),f(b)\rangle+\langle f(a),y\rangle
 =\mathcal B_{\phi}
    \big((\mathrm{id}_{\frak g}+f)(x+a),
    (\mathrm{id}_{\frak g}+f)(y+b)\big).
    \end{equation*}
    Hence
    $\big((\frak g\oplus\frak g^*,[-,-]'_d,\mathcal B_{\phi}),
    (\frak g,[-,-]_{\frak g}),
    (\frak g^*,[-,-]_{\frak g^*})\big)$
    is a standard generic Manin triple of Lie algebras,
    and the two generic Manin triples are isomorphic through
    $\mathrm{id}_{\frak g}+f:\frak g\oplus\frak h
    \rightarrow\frak g\oplus\frak g^*$.
\end{proof}

(Standard) Manin triples of Lie algebras with respect to the nondegenerate symmetric invariant bilinear forms, symplectic forms, or nondegenerate commutative $2$-cocycles are all typical
 examples of (standard) generic Manin triples of Lie algebras (Definitions \ref{2388}, \ref{2388-} and \ref{2388--}).

\begin{lem}\label{1452}
    Let $\big((\frak g\oplus\frak g^*,[-,-]_d,\mathcal B_{\phi}),
    (\frak g,[-,-]_{\frak g}),
    (\frak g^*,[-,-]_{\frak g^*})\big)$
    be a standard generic Manin triple of Lie algebras.
    Let $P_{\frak g}:\frak g\oplus\frak g^*\rightarrow\frak g$
    and $P_{\frak g^*}:\frak g\oplus\frak g^*\rightarrow\frak g^*$
    be the projections given by
    \begin{equation*}\label{1453}
        P_{\frak g}(x+\xi )=x,\quad
        P_{\frak g^*}(x+\xi )=\xi ,
        \quad \forall x\in\frak g,\ \xi \in\frak g^*.
    \end{equation*}
    Let $\circ_{\frak g}:\frak g\otimes\frak g\rightarrow\frak g$
    and $\circ_{\frak g^*}:\frak g^*\otimes\frak g^*\rightarrow\frak g^*$
    be multiplications given through
    \begin{equation}\label{1458}
        L^*_{\circ_{\frak g}}(x)\xi
        =P_{\frak g^*}([x,\xi ]_d),\quad
        L^*_{\circ_{\frak g^*}}(\xi )x
        =P_{\frak g}([\xi ,x]_d),
        \quad \forall x\in\frak g,\ \xi \in\frak g^*,
    \end{equation}
    or equivalently,
    \begin{align*}
        \langle x\circ_{\frak g}y,\xi \rangle
        =-\langle y,P_{\frak g^*}([x,\xi ]_d)\rangle,\;
        \langle \xi \circ_{\frak g^*}\eta ,x\rangle
        =-\langle \eta ,P_{\frak g}([\xi ,x]_d)\rangle,\;\forall x,y\in\frak g,\xi ,\eta \in\frak g^*.
    \end{align*}
    Then the Lie algebra structure
    $(\frak g\oplus\frak g^*,[-,-]_d)$ can be explicitly written as
    \begin{equation}\label{2662}
        [x+\xi ,y+\eta ]_d
        =[x,y]_{\frak g}
        +L^*_{\circ_{\frak g^*}}(\xi )y
        -L^*_{\circ_{\frak g^*}}(\eta )x+[\xi ,\eta ]_{\frak g^*}
        +L^*_{\circ_{\frak g}}(x)\eta
        -L^*_{\circ_{\frak g}}(y)\xi ,
    \end{equation}
    for all $x,y\in\frak g$ and $\xi ,\eta \in\frak g^*$.
    That is, $(\frak g,[-,-]_{\frak g},\circ_{\frak g})$ and
    $(\frak g^*,[-,-]_{\frak g^*},\circ_{\frak g^*})$
    are self-module Lie algebras, and
    $\big((\frak g,[-,-]_{\frak g}),
    (\frak g^*,[-,-]_{\frak g^*}),
    L^*_{\circ_{\frak g}},L^*_{\circ_{\frak g^*}}\big)$
    is a matched pair of Lie algebras.
    Conversely, suppose that
    $(\frak g,[-,-]_{\frak g},\circ_{\frak g})$ and
    $(\frak g^*,[-,-]_{\frak g^*},\circ_{\frak g^*})$
    are self-module Lie algebras.
    If $\big((\frak g,[-,-]_{\frak g}),
    (\frak g^*$,
    $[-,-]_{\frak g^*}),
    L^*_{\circ_{\frak g}},L^*_{\circ_{\frak g^*}}\big)$
    is a matched pair of Lie algebras, then
    $\big((\frak g\oplus\frak g^*,[-,-]_d,\mathcal B_{\phi}),
    (\frak g,[-,-]_{\frak g}),
    (\frak g^*$,
    $[-,-]_{\frak g^*})\big)$
    with $[-,-]_d$ given by \eqref{2662}
    is a standard generic Manin triple of Lie algebras for any linear
    isomorphism $\phi $ on $\frak g$.
\end{lem}

\begin{proof}
    Suppose that
    $\big((\frak g\oplus\frak g^*,[-,-]_d,\mathcal B_{\phi}),
    (\frak g,[-,-]_{\frak g}),
    (\frak g^*,[-,-]_{\frak g^*})\big)$
    is a standard generic Manin triple of Lie algebras.
    For all $x,y\in\frak g$ and $\xi ,\eta \in\frak g^*$, we have
    \begin{align*}
        [x+\xi ,y+\eta ]_d
        &=[x,y]_{\frak g}+[x,\eta ]_d
        +[\xi ,y]_d+[\xi ,\eta ]_{\frak g^*}\\
        &\overset{\eqref{1458}}{=}
        [x,y]_{\frak g}
        +L^*_{\circ_{\frak g^*}}(\xi )y
        -L^*_{\circ_{\frak g^*}}(\eta )x+[\xi ,\eta ]_{\frak g^*}
        +L^*_{\circ_{\frak g}}(x)\eta
        -L^*_{\circ_{\frak g}}(y)\xi .
    \end{align*}
    Then by Definition \ref{1493},
    $\big((\frak g,[-,-]_{\frak g}),
    (\frak g^*,[-,-]_{\frak g^*}),
    L^*_{\circ_{\frak g}},L^*_{\circ_{\frak g^*}}\big)$
    is a matched pair of Lie algebras.
    Conversely, suppose that the above quadruple is a matched
    pair of Lie algebras. Then $[-,-]_d$ given by \eqref{2662}
    defines a Lie algebra structure on $\frak g\oplus\frak g^*$
    containing $\frak g$ and $\frak g^*$ as Lie subalgebras.
    For any linear isomorphism $\phi$ on $\frak g$,
    $\mathcal B_{\phi}$ is nondegenerate, and both summands
    are isotropic with respect to it.
    Hence
    $\big((\frak g\oplus\frak g^*,[-,-]_d,\mathcal B_{\phi}),
    (\frak g,[-,-]_{\frak g}),
    (\frak g^*,[-,-]_{\frak g^*})\big)$
    is a standard generic Manin triple of Lie algebras.
\end{proof}

\begin{pro}\label{1503}
    Suppose that  $\big((\frak g\oplus\frak g^*,[-,-]_d,\mathcal B_{\phi}),
    (\frak g,[-,-]_{\frak g}),
    (\frak g^*,[-,-]_{\frak g^*})\big)$
    is a standard
    generic Manin triple of Lie algebras. Let
    $(\frak g,[-,-]_{\frak g},\circ_{\frak g})$ and
    $(\frak g^*,[-,-]_{\frak g^*},\circ_{\frak g^*})$
    be self-module Lie algebras given in Lemma \ref{1452}
    that make
    $\big((\frak g,[-,-]_{\frak g}),
    (\frak g^*,[-,-]_{\frak g^*}),
    L^*_{\circ_{\frak g}},L^*_{\circ_{\frak g^*}}\big)$ into a matched pair of Lie algebras.
    Let $(\frak g\oplus\frak g^*,[-,-]_d,\bullet_{\phi})$
    be the self-module Lie algebra induced from
    $(\frak g\oplus\frak g^*,[-,-]_d,\mathcal B_{\phi})$.
    Then the multiplication $\bullet_{\phi}$ is explicitly given by
    \begin{align}
        (x+\xi )\bullet_{\phi}(y+\eta )
        &=x\circ_\phi y
        +\mathrm{ad}^*_{\frak g^*}[\phi^{-1},1,\phi](\xi )y
        +R^*_{\circ_{\frak g^*}}[\phi^{-1},1,1](\eta )x
        \nonumber\\
        &\quad
        +\xi \circ_{\frak g^*}\eta
        +\mathrm{ad}^*_{\frak g}(x)\eta
        +R^*_{\circ_{\frak g}}[1,\phi,1](y)\xi ,
        \label{1125-}
    \end{align}
    for all $x,y\in\frak g$ and $\xi ,\eta \in\frak g^*$,
    where $(\frak g,[-,-]_{\frak g},\circ_\phi)$
    is the $\phi$-conjugate self-module Lie algebra of
    $(\frak g,[-,-]_{\frak g},\circ_{\frak g})$.
    Therefore, the collection
    \begin{equation}\label{1529}
        \big(
        (\frak g,[-,-]_{\frak g},\circ_\phi),
        (\frak g^*,[-,-]_{\frak g^*},\circ_{\frak g^*}),
     L^*_{\circ_{\frak g}},
        \mathrm{ad}^*_{\frak g},
        R^*_{\circ_{\frak g}}[1,\phi,1],
         L^*_{\circ_{\frak g^*}},
        \mathrm{ad}^*_{\frak g^*}[\phi^{-1},1,\phi],
        R^*_{\circ_{\frak g^*}}[\phi^{-1},1,1]
        \big)
    \end{equation} is a matched pair of self-module Lie algebras.
\end{pro}

\begin{proof}
    For all $x,y,z\in\frak g$ and $\xi \in\frak g^*$, we have
    \begin{align*}
        \mathcal B_{\phi}(x\bullet_{\phi}y,z+\xi )
        &=-\mathcal B_{\phi}(y,[x,z+\xi ]_d) =-\langle\phi(y),L^*_{\circ_{\frak g}}(x)\xi \rangle
        =\langle x\circ_{\frak g}\phi(y),\xi \rangle\\
        &=\mathcal B_{\phi}(x\circ_\phi y,z+\xi ).
    \end{align*}
    By the nondegeneracy of $\mathcal B_{\phi}$,
    we have $x\bullet_{\phi}y=x\circ_\phi y$.
    Similarly we have
    $\xi \bullet_{\phi}\eta =\xi \circ_{\frak g^*}\eta $
    for all $\xi ,\eta \in\frak g^*$.
    Moreover, for all $x,y\in\frak g$ and
    $\xi ,\eta \in\frak g^*$, we have
    \begin{align*}
        \mathcal B_{\phi}(x\bullet_{\phi}\xi ,y+\eta )
        &=-\mathcal B_{\phi}(\xi ,[x,y+\eta ]_d) =\langle\mathrm{ad}^*_{\frak g}(x)\xi ,y\rangle
        -\langle \eta \circ_{\frak g^*}\xi ,x\rangle\\
        &=\mathcal B_{\phi}
        \big(
        \mathrm{ad}^*_{\frak g}(x)\xi
        +\phi^{-1}R^*_{\circ_{\frak g^*}}(\xi )x,
        y+\eta
        \big).
    \end{align*}
    Hence
    $x\bullet_{\phi}\xi
    =\mathrm{ad}^*_{\frak g}(x)\xi
    +\phi^{-1}R^*_{\circ_{\frak g^*}}(\xi )x$.
    Similarly,
    $\xi \bullet_{\phi}x
    =\phi^{-1}\mathrm{ad}^*_{\frak g^*}(\xi )\phi(x)
    +R^*_{\circ_{\frak g}}(\phi(x))\xi $.
    Therefore,  \eqref{1125-} holds. Then
    by Definition \ref{defi:mpSR-Lie}, the collection
    \eqref{1529} is a matched pair of self-module Lie algebras.
\end{proof}

\begin{pro}\label{1559}
    Let $(\frak g,[-,-]_{\frak g}),(\frak g^*,[-,-]_{\frak g^*})$ be Lie algebras and $\phi$ be an arbitrary linear isomorphism on $\frak g$.
    Then the following statements are equivalent.
    \begin{enumerate}
        \item\label{1563}
        There is a Lie algebra structure $(\frak g\oplus\frak g^*,[-,-]_d)$ on $\frak g\oplus\frak g^{*}$ such that
        $\big((\frak g\oplus\frak g^*,[-,-]_d$,
        $\mathcal B_{\phi}),
        (\frak g,[-,-]_{\frak g}),(\frak g^*,[-,-]_{\frak g^*})\big)$
        is a standard generic Manin triple of Lie algebras.

        \item\label{1564}
        There are self-module Lie algebras
        $(\frak g,[-,-]_{\frak g},\circ_{\frak g}),(\frak g^*,[-,-]_{\frak g^*},\circ_{\frak g^*})$
        such that
        $\big((\frak g,[-,-]_{\frak g})$,
        $(\frak g^*,[-,-]_{\frak g^*}),
        L^*_{\circ_{\frak g}},L^*_{\circ_{\frak g^*}}\big)$
        is a matched pair of Lie algebras.

        \item\label{1565}
        There are self-module Lie algebras
        $(\frak g,[-,-]_{\frak g},\circ_{\frak g}),(\frak g^*,[-,-]_{\frak g^*},\circ_{\frak g^*})$
        such that the collection \eqref{1529} is a matched pair of self-module Lie algebras.
    \end{enumerate}
\end{pro}

\begin{proof}
    By Lemma \ref{1452} and Proposition \ref{1503}, we have
    \eqref{1563}$\Longrightarrow$\eqref{1565}.
    By Proposition \ref{pro:mat}, we have
    \eqref{1565}$\Longrightarrow$\eqref{1564}.
    By Lemma \ref{1452} again, we have
    \eqref{1564}$\Longrightarrow$\eqref{1563}.
\end{proof}

\subsection{Self-module Lie bialgebras}

\begin{defi}
    Let $\frak g$ be a vector space with comultiplications $\Delta,\theta :\frak g\rightarrow \frak g\otimes\frak g$.
    Suppose that $(\frak g,\Delta)$ is a {\bf Lie coalgebra}, that is, the following equations hold:
    \begin{eqnarray}
        &&(\Delta+\tau\Delta)(x)=0,\label{eq:Lie coalg1}\\
        &&\sigma(\mathrm{id}_{\frak g}\otimes \Delta)\Delta(x)=0,\;\forall x\in\frak g,\label{eq:Lie coalg2}
    \end{eqnarray}
where $\sigma:\frak g\otimes\frak g\otimes\frak g\rightarrow\frak g\otimes\frak g\otimes\frak g$ is a linear map given by
\begin{equation*}
     \sigma(x\otimes y\otimes z)=x\otimes y\otimes z+y\otimes z\otimes x+z\otimes x\otimes y,\;\forall x,y,z\in \frak g.
\end{equation*}
    If the following equation holds:
    \begin{eqnarray}
        (\Delta\otimes\mathrm{id}_{\frak g})\theta(x)= (\mathrm{id}_{\frak g}\otimes\theta)\theta(x)-(\tau\otimes\mathrm{id}_{\frak g})(\mathrm{id}_{\frak g}\otimes\theta)\theta(x),\;\forall x\in\frak g,\label{eq:NS-coalg2}
    \end{eqnarray}
    then we say $(\frak g,\Delta,\theta)$ is a {\bf self-module Lie coalgebra}.
\end{defi}

\begin{pro}\label{pro:2368}
    Let $\frak g$ be a vector space with comultiplications $\Delta,\theta:\frak g\rightarrow \frak g\otimes \frak g$.
    Let $[-,-]_{\frak g^{*}} ,\circ_{\frak g^{*}} :\frak g^{*}\otimes \frak g^{*}\rightarrow \frak g^{*}$
    be the linear duals of $\Delta$ and $\theta$ respectively, that is,
    \begin{equation}\label{eq:2249}
     \langle [\xi ,\eta ]_{\frak g^{*}},x\rangle =\langle \xi \otimes \eta ,\Delta(x)\rangle, \;\langle \xi \circ_{\frak g^{*}} \eta ,x\rangle =\langle \xi \otimes \eta ,\theta(x)\rangle,\;\forall x\in \frak g, \xi ,\eta \in \frak g^{*}.
    \end{equation}
    Then $(\frak g,\Delta,\theta)$ is a self-module Lie coalgebra if and only if $(\frak g^{*},[-,-]_{\frak g^{*}},\circ_{\frak g^{*}} )$ is a self-module Lie algebra.
\end{pro}
\begin{proof}
    The conclusion follows by dualizing the defining identities.
\end{proof}

\begin{defi}
    Let $(\frak g,[-,-]_{\frak g},\circ_{\frak g} )$ be a self-module Lie algebra and $(\frak g,\Delta,\theta)$ be a self-module Lie coalgebra.
    If the following equations hold for all $x,y\in\frak g$:
    \begin{align}
        \Delta(x\circ_{\frak g}y)&=(\mathrm{id}_{\frak g}^{\otimes 2}-\tau)\Big(\big( L_{\circ_{\frak g}}(x)\otimes\mathrm{id}_{\frak g} \big)\Delta(y)+ \big(\mathrm{id}_{\frak g}\otimes R_{\circ_{\frak g}}(y)\big) \theta(x)\Big),\label{eq:NS bialg1}\\
        \theta([x,y]_{\frak g})&=\big( L_{\circ_{\frak g}}(x)\otimes\mathrm{id}_{\frak g}+\mathrm{id}_{\frak g}\otimes\mathrm{ad}_{\frak g}(x)  \big)\theta(y)
        -\big( L_{\circ_{\frak g}}(y)
        \otimes\mathrm{id}_{\frak g}+\mathrm{id}_{\frak g}\otimes\mathrm{ad}_{\frak g}(y)  \big)\theta(x)
        , \label{eq:NS bialg2}
    \end{align}
    then we say $(\frak g,[-,-]_{\frak g},\circ_{\frak g} ,\Delta,\theta)$ is a {\bf self-module Lie bialgebra}.
\end{defi}

Lie bialgebras, pre-Lie bialgebras and anti-pre-Lie bialgebras are
typical examples of self-module Lie bialgebras (Propositions \ref{2920}, \ref{2646} and \ref{2932}). In
particular, a self-module Lie bialgebra is not necessarily a Lie
bialgebra.

Recall that a \textbf{$1$-cocycle} of a Lie algebra
$(\mathfrak{g},[-,-]_{\frak g})$ associated to a representation $(\rho,V)$ is a linear map
$f:\mathfrak{g}\rightarrow V$ such that the following equation holds:
\begin{equation*}
    f([x,y]_{\frak g})=\rho(x)f(y)-\rho(y)f(x),\;\;\forall x,y\in\mathfrak{g}.
\end{equation*}
Given two representations $(\rho,V)$ and $(\mu,W)$ of a Lie algebra
$(\mathfrak{g},[-,-]_{\frak g})$, it is known that
$(\rho\otimes\mathrm{id}_{W}+\mathrm{id}_{V}\otimes\mu, V\otimes W)$ is
also a representation of $(\mathfrak{g},[-,-]_{\frak g})$, where
\begin{equation*}
    (\rho\otimes\mathrm{id}_{W}+\mathrm{id}_{V}\otimes\mu)(x)(v\otimes w)
    =\rho(x)v\otimes w+v\otimes\mu(x)w, \;\;\forall x\in\mathfrak{g}, v\in V,w\in W.
\end{equation*}

\begin{pro}\label{pro:equ2}
    Let $(\frak g,[-,-]_{\frak g} )$ and $(\frak g^{*},[-,-]_{\frak g^{*}} )$ be Lie algebras and $\circ_{\frak g}:\frak g\otimes\frak g\rightarrow \frak g, \circ_{\frak g^{*}}:\frak g^{*}\otimes\frak g^{*}\rightarrow\frak g^{*}$ be multiplications.
    Let
    $\Delta,\theta :\frak g\rightarrow\frak g\otimes\frak g$ be comultiplications given by \eqref{eq:2249}.
    Then the following statements are equivalent:
    \begin{enumerate}
        \item\label{2787} $\big( (\frak g,[-,-]_{\frak g}),(\frak g^{*},[-,-]_{\frak g^{*}}),L^{*}_{\circ_{\frak g}}, L^{*}_{\circ_{\frak g^{*}}} \big)$ is a matched pair of Lie algebras.
        \item\label{2789} $(L_{\circ_{\frak g}},\frak g)$ and $(L_{\circ_{\frak g^{*}}},\frak g^{*})$ are representations of $(\frak g,[-,-]_{\frak g})$ and $(\frak g^{*},[-,-]_{\frak g^{*}})$ respectively, $\theta$ is a $1$-cocycle of $(\frak g,[-,-]_{\frak g})$ associated to $(L_{\circ_{\frak g}}\otimes\mathrm{id}_{\frak g}+\mathrm{id}_{\frak g}\otimes\mathrm{ad}_{\frak g},\frak g\otimes\frak g)$ and $\vartheta$ is a $1$-cocycle of $(\frak g^{*},[-,-]_{\frak g^{*}})$ associated to $(L_{\circ_{\frak g^{*}}}\otimes\mathrm{id}_{\frak g ^* }+\mathrm{id}_{ \frak g^* }\otimes\mathrm{ad}_{\frak g^{*}},\frak g^{*}\otimes\frak g^{*})$, where $\vartheta:\frak g^{*}\rightarrow\frak g^{*}\otimes\frak g^{*}$ is the linear dual of $\circ_{\frak g}$.
        \item\label{2793} $(\frak g,[-,-]_{\frak g},\circ_{\frak g},\Delta,\theta)$ is a self-module Lie bialgebra.
    \end{enumerate}
\end{pro}
\begin{proof}
    \eqref{2787}$\Longleftrightarrow$\eqref{2789}.  Suppose that $\big( (\frak g,[-,-]_{\frak g}),(\frak g^{*},[-,-]_{\frak g^{*}}),L^{*}_{\circ_{\frak g}}, L^{*}_{\circ_{\frak g^{*}}} \big)$ is a matched pair of Lie algebras. Then $(L^{*}_{\circ_{\frak g }},\frak g^{*})$
    as well as $(L_{\circ_{\frak g}},\frak g)$ is a representation of $(\frak g,[-,-]_{\frak g})$,
    $(L^{*}_{\circ_{\frak g^{*}}},\frak g )$
    as well as $(L_{\circ_{\frak g^{*}}},\frak g^{*})$ is a representation of $(\frak g^{*},[-,-]_{\frak g^{*}})$, and the following equations hold for all $x,y\in\frak g, \xi ,\eta \in\frak g^{*}$:
        \begin{align}
        L_{\circ_{\frak g^{*}}}^{*}(\xi )[x,y]_{\frak g}-[L_{\circ_{\frak g^{*}}}^{*}(\xi )x,y]_{\frak g}
        -[x,L_{\circ_{\frak g^{*}}}^{*}(\xi )y]_{\frak g}+L_{\circ_{\frak g^{*}}}^{*}\big( L_{\circ_{\frak g}}^{*}(x)\xi  \big)y-L_{\circ_{\frak g^{*}}}^{*}\big( L_{\circ_{\frak g}}^{*}(y)\xi  \big)x&=0, \label{2813}\\
        L_{\circ_{\frak g}}^{*}(x)[\xi ,\eta ]_{\frak g^{*}}-[L_{\circ_{\frak g}}^{*}(x)\xi ,\eta ]_{\frak g^{*}}
        -[\xi , L_{\circ_{\frak g}}^{*}(x)\eta ]_{\frak g^{*}}+L_{\circ_{\frak g}}^{*}\big( L_{\circ_{\frak g^{*}}}^{*}(\xi )x \big)\eta -L_{\circ_{\frak g}}^{*}\big( L_{\circ_{\frak g^{*}}}^{*}(\eta )x \big)\xi &=0. \label{2815}
    \end{align}
For all $x,y\in\frak g, \xi , \eta \in\frak g^{*}$, we have
    \begin{align*}
        \langle L_{\circ_{\frak g^{*}}}^{*}(\xi ) [x,y]_{\frak g} , \eta \rangle&= -\langle [x,y]_{\frak g}, \xi  \circ_{\frak g^{*}}\eta \rangle=-\langle \theta([x,y]_{\frak g}), \xi \otimes \eta \rangle,\\
        -\langle [L_{\circ_{\frak g^{*}}}^{*}(\xi )x,y]_{\frak g},\eta \rangle&= -\langle L_{\circ_{\frak g^{*}}}^{*}(\xi )x,\mathrm{ad}^{*}_{\frak g}(y)\eta \rangle=\langle x, \xi  \circ_{\frak g^{*}}\mathrm{ad}^{*}_{\frak g}(y)\eta \rangle=\langle \theta(x), \xi \otimes\mathrm{ad}^{*}_{\frak g}(y)\eta \rangle\\
        &= -\langle \big(\mathrm{id}_{\frak g}\otimes\mathrm{ad}_{\frak g}(y)\big)\theta(x), \xi \otimes \eta \rangle,\\
        \langle L_{\circ_{\frak g^{*}}}^{*}\big(L_{\circ_{\frak g}}^{*}(x)\xi \big) y, \eta \rangle
        &=-\langle \theta(y),L_{\circ_{\frak g}}^{*}(x)\xi \otimes \eta \rangle = \langle \big(L_{\circ_{\frak g}}(x)\otimes\mathrm{id}_{\frak g} \big)\theta(y), \xi \otimes \eta \rangle,
    \end{align*}
    and similarly
    \begin{align*}
        -\langle [x,L_{\circ_{\frak g^{*}}}^{*}(\xi )y]_{\frak g}, \eta \rangle&= \langle \big(\mathrm{id}_{\frak g}\otimes\mathrm{ad}_{\frak g}(x)\big)\theta(y), \xi \otimes \eta  \rangle, \\
        -\langle L_{\circ_{\frak g^{*}}}^{*}\big(L_{\circ_{\frak g }}^{*}(y)\xi \big) x, \eta \rangle&= -\langle \big( L_{\circ_{\frak g}}(y)\otimes\mathrm{id}_{\frak g} \big)\theta(x), \xi \otimes \eta \rangle.
    \end{align*}
    Thus by \eqref{2813}, $\theta$ is a $1$-cocycle of $(\frak g,[-,-]_{\frak g})$ associated to $(L_{\circ_{\frak g}}\otimes\mathrm{id}_{\frak g}+\mathrm{id}_{\frak g}\otimes\mathrm{ad}_{\frak g},\frak g\otimes\frak g)$. Similarly by \eqref{2815}, $\vartheta$ is a $1$-cocycle of $(\frak g^{*},[-,-]_{\frak g^{*}})$ associated to $(L_{\circ_{\frak g^{*}}}\otimes\mathrm{id}_{\frak g^{*}}+\mathrm{id}_{\frak g^{*}}\otimes\mathrm{ad}_{\frak g^{*}},\frak g^{*}\otimes\frak g^{*})$.
    The converse implication follows similarly.

    \eqref{2789}$\Longleftrightarrow$\eqref{2793}. Suppose that the assumption \eqref{2789} holds.
    It is clear that $(\frak g,[-,-]_{\frak g},\circ_{\frak g})$
    and $(\frak g^{*},[-,-]_{\frak g^{*}}$,
    $\circ_{\frak g^{*} })$ are self-module Lie algebras, and hence $(\frak g,\Delta,\theta)$ is a self-module Lie coalgebra.
    Since $\theta$ is a $1$-cocycle of $(\frak g,[-,-]_{\frak g})$ associated to $(L_{\circ_{\frak g}}\otimes\mathrm{id}_{\frak g}+\mathrm{id}_{\frak g}\otimes\mathrm{ad}_{\frak g},\frak g\otimes\frak g)$, \eqref{eq:NS bialg2} holds clearly. For all $x,y\in\frak g, \xi ,\eta \in\frak g^{*}$, we have
\begin{align*}
	 \langle \vartheta([\xi ,\eta ]_{\frak g^{*}}), x\otimes y\rangle&= \langle [\xi ,\eta ]_{\frak g^{*}}, x\circ_{\frak g}y\rangle\\
	 &= \langle \xi \otimes \eta ,\Delta ( x\circ_{\frak g}y)\rangle,\\
	- \langle \big( L_{\circ_{\frak g^{*}}}(\eta )\otimes\mathrm{id}_{\frak g^{*}}
	+\mathrm{id}_{\frak g^{*}}\otimes \mathrm{ad}_{\frak g^{*}}(\eta ) \big)\vartheta(\xi ),x\otimes y\rangle 
	&= \langle \xi ,x\circ_{\frak g}\mathrm{ad}^{*}_{\frak g^{*}}(\eta )y+L_{\circ_{\frak g^{*}}}^{*}(\eta )x \circ_{\frak g}y\rangle\\
	&= -\langle [L_{\circ_{\frak g}}^{*}(x)\xi ,\eta ]_{\frak g^{*}},y\rangle-\langle R_{\circ_{\frak g}}^{*}(y)\xi ,L_{\circ_{\frak g^{*}}}^{*}(\eta )x \rangle\\
	&= \langle \xi \otimes \eta , \big( L_{\circ_{\frak g}}(x)\otimes\mathrm{id}_{\frak g} \big)\Delta(y)-\tau \big( \mathrm{id}_{\frak g}\otimes R_{\circ_{\frak g}}(y) \big)\theta(x)\rangle,
\end{align*}
    and similarly
    \begin{equation*}
        \langle \big( L_{\circ_{\frak g^{*}}}(\xi )\otimes\mathrm{id}_{\frak g^{*}}
        +\mathrm{id}_{\frak g^{*}}\otimes \mathrm{ad}_{\frak g^{*}}(\xi ) \big)\vartheta(\eta ), x\otimes y\rangle
        =\langle \xi \otimes \eta , \big( \mathrm{id}_{\frak g}\otimes R_{\circ_{\frak g}}(y) \big)\theta(x)-\tau\big( L_{\circ_{\frak g}}(x)\otimes\mathrm{id}_{\frak g} \big)\Delta(y) \rangle.
    \end{equation*}
Since $\vartheta$ is a $1$-cocycle of $(\frak g^{*},[-,-]_{\frak g^{*}})$ associated to $(L_{\circ_{\frak g^{*}}}\otimes\mathrm{id}_{\frak g^{*}}+\mathrm{id}_{\frak g^{*}}\otimes\mathrm{ad}_{\frak g^{*}},\frak g^{*}\otimes\frak g^{*})$, \eqref{eq:NS bialg1} holds. Hence we obtain \eqref{2793}.
    The converse implication follows similarly.
\end{proof}

Combining Propositions \ref{1559} and \ref{pro:equ2}, we obtain
the following result.

\begin{thm}\label{cor:2532}
    Let $(\frak g,[-,-]_{\frak g}),(\frak g^*,[-,-]_{\frak g^*})$ be Lie algebras and $\phi$ be an arbitrary linear isomorphism on $\frak g$.
    Then the following statements are equivalent:
    \begin{enumerate}
        \item\label{1563-}
    There is a Lie algebra structure
    $(\frak g\oplus\frak g^*,[-,-]_d)$ such that
    $\big((\frak g\oplus\frak g^*,[-,-]_d,\mathcal B_{\phi}),
    (\frak g,[-,-]_{\frak g})$,
    $(\frak g^*,[-,-]_{\frak g^*})\big)$
    is a standard generic Manin triple of Lie algebras.

        \item\label{1564-}
        There are self-module Lie algebras
        $(\frak g,[-,-]_{\frak g},\circ_{\frak g}),(\frak g^*,[-,-]_{\frak g^*},\circ_{\frak g^*})$
        such that
        $\big((\frak g,[-,-]_{\frak g})$,
        $(\frak g^*,[-,-]_{\frak g^*}),
        L^*_{\circ_{\frak g}},L^*_{\circ_{\frak g^*}}\big)$
        is a matched pair of Lie algebras.

        \item\label{1566-}
        There are self-module Lie algebras
        $(\frak g,[-,-]_{\frak g},\circ_{\frak g}),(\frak g^*,[-,-]_{\frak g^*},\circ_{\frak g^*})$
        such that the collection \eqref{1529}
        is a matched pair of self-module Lie algebras.

        \item\label{1567-}
        There are self-module Lie algebras
        $(\frak g,[-,-]_{\frak g},\circ_{\frak g}),(\frak g^*,[-,-]_{\frak g^*},\circ_{\frak g^*})$
        such that
        $(\frak g,[-,-]_{\frak g}$,
        $\circ_{\frak g},\Delta,\theta)$
        is a self-module Lie bialgebra, where
        $\Delta,\theta:\frak g\rightarrow\frak g\otimes\frak g$
        are given by \eqref{eq:2249}.
    \end{enumerate}
\end{thm}

\begin{rmk}
    One may find the equivalence \eqref{1566-}$\Longleftrightarrow$\eqref{1567-} in Theorem \ref{cor:2532} to some
    extent counter-intuitive, since the matched pair of self-module
    Lie algebras in \eqref{1566-} involves the parameter $\phi$,
    whereas the corresponding self-module Lie bialgebra in
    \eqref{1567-} involves no such parameter.
    To clarify this, we present the following explicit computational proof of \eqref{1566-}$\Longleftrightarrow$\eqref{1567-}.
    For all $x,y\in\frak g$ and $\xi ,\eta \in\frak g^*$, we have
    \begin{align*}
        \langle R^*_{\circ_{\frak g^*}}[\phi^{-1},1,1](\xi )[x,y]_{\frak g},\eta \rangle
        &=\langle R^*_{\circ_{\frak g^*}}(\xi )[x,y]_{\frak g},(\phi^{-1})^t(\eta )\rangle =-\langle [x,y]_{\frak g},(\phi^{-1})^t(\eta )\circ_{\frak g^*}\xi \rangle\\
        &=-\langle (\phi^{-1}\otimes\mathrm{id}_{\frak g})\theta([x,y]_{\frak g}),\eta \otimes \xi \rangle.
    \end{align*}
    Similarly, we have
    \begin{align*}
        -\langle x\circ_\phi\big(R^*_{\circ_{\frak g^*}}[\phi^{-1},1,1](\xi )y\big),\eta \rangle
        &=\langle (\phi^{-1}\otimes\mathrm{id}_{\frak g})\big(L_{\circ_{\frak g}}(x)\otimes\mathrm{id}_{\frak g}\big)\theta(y),\eta \otimes \xi \rangle,\\
        -\langle R^*_{\circ_{\frak g^*}}[\phi^{-1},1,1]\big(\mathrm{ad}^*_{\frak g}(y)\xi \big)x,\eta \rangle
        &=-\langle (\phi^{-1}\otimes\mathrm{id}_{\frak g})\big(\mathrm{id}_{\frak g}\otimes\mathrm{ad}_{\frak g}(y)\big)\theta(x),\eta \otimes \xi \rangle,\\
        \langle y\circ_\phi\big(R^*_{\circ_{\frak g^*}}[\phi^{-1},1,1](\xi )x\big),\eta \rangle
        &=-\langle (\phi^{-1}\otimes\mathrm{id}_{\frak g})\big(L_{\circ_{\frak g}}(y)\otimes\mathrm{id}_{\frak g}\big)\theta(x),\eta \otimes \xi \rangle,\\
        \langle R^*_{\circ_{\frak g^*}}[\phi^{-1},1,1]\big(\mathrm{ad}^*_{\frak g}(x)\xi \big)y,\eta \rangle
        &=\langle (\phi^{-1}\otimes\mathrm{id}_{\frak g})\big(\mathrm{id}_{\frak g}\otimes\mathrm{ad}_{\frak g}(x)\big)\theta(y),\eta \otimes \xi \rangle.
    \end{align*}
    Since the common factor $\phi^{-1}\otimes\mathrm{id}_{\frak g}$ is invertible, the collection \eqref{1529} satisfies \eqref{12641}
    if and only if \eqref{eq:NS bialg2} holds.
    Similarly, regarding the same collection \eqref{1529},  we have \eqref{12642}$\Longleftrightarrow$\eqref{eq:NS bialg1}$\Longleftrightarrow$\eqref{12643} and \eqref{12644}$\Longleftrightarrow$\eqref{eq:NS bialg2}.
\end{rmk}

\vspace{-.3cm}

\section{The self-module Lie Yang-Baxter equation and $\mathcal{O}$-operators}\label{sec:4}
We introduce the self-module Lie Yang-Baxter equation (SLYBE), whose solutions yield self-module Lie bialgebras.
Such a bialgebra is triangular if the solution satisfies the $\phi$-condition for a linear isomorphism $\phi$ on $\frak g$.
We characterize  solutions of the SLYBE as $\mathcal{O}$-operators of the underlying Lie algebra and, under the above condition, of the $\phi$-conjugate self-module Lie algebra.
Finally, we use self-bimodule pre-Lie algebras to construct $\mathcal{O}$-operators of self-module Lie algebras and hence SLYBE solutions in the corresponding semi-direct product self-module Lie algebras.

\subsection{The self-module Lie Yang-Baxter equation}

\begin{defi}
    Let $(\frak g,[-,-]_{\frak g},\circ_{\frak g})$ be a self-module Lie algebra and $r=\sum\limits_i a_i\otimes b_i\in\frak g\otimes\frak g$.
    Set ${\bf SL}(r)\in\frak g\otimes\frak g\otimes\frak g$ by
    \begin{equation*}
        {\bf SL}(r)=r_{12}\circ_{\frak g}r_{23}-r_{21}\circ_{\frak g}r_{13}+[r_{13},r_{23}]_{\frak g},
    \end{equation*}
    where
    \begin{equation*}
        r_{12}\circ_{\frak g}r_{23}=\sum_{i,j}a_i\otimes b_i\circ_{\frak g}a_j\otimes b_j,\;
        r_{21}\circ_{\frak g}r_{13}=\sum_{i,j}b_i\circ_{\frak g}a_j\otimes a_i\otimes b_j,\;
        [r_{13},r_{23}]_{\frak g}=\sum_{i,j}a_i\otimes a_j\otimes[b_i,b_j]_{\frak g}.
    \end{equation*}
    We say $r$ is a solution of the {\bf self-module Lie Yang-Baxter equation} (or {\bf SLYBE} in short) in $(\frak g,[-,-]_{\frak g},\circ_{\frak g})$ if ${\bf SL}(r)=0$.
\end{defi}

\begin{pro}\label{pro:4.2}
    Let $(\frak g,[-,-]_{\frak g},\circ_{\frak g})$ be a self-module Lie algebra and $r=\sum\limits_i a_i\otimes b_i\in\frak g\otimes\frak g$.
    Let $\Delta,\theta:\frak g\rightarrow\frak g\otimes\frak g$ be comultiplications given by
        \begin{align}
            \Delta(x)&=\Delta_r(x)=(\tau-\mathrm{id}_{\frak g}^{\otimes 2})\big(\mathrm{id}_{\frak g}\otimes R_{\circ_{\frak g}}(x)\big)r
            =\sum_i b_i\circ_{\frak g}x\otimes a_i-a_i\otimes b_i\circ_{\frak g}x,\label{eq:NS com1}\\
            \theta(x)&=\theta_r(x)=\big(L_{\circ_{\frak g}}(x)\otimes\mathrm{id}_{\frak g}+\mathrm{id}_{\frak g}\otimes\mathrm{ad}_{\frak g}(x)\big)r
            =\sum_i x\circ_{\frak g}a_i\otimes b_i+a_i\otimes[x,b_i]_{\frak g},\;\forall x\in\frak g.\label{eq:NS com2}
        \end{align}
    \begin{enumerate}
        \item\label{1530} $(\frak g,\Delta)$ is a Lie coalgebra if and only if the following equation holds:
        \begin{eqnarray}\label{eq:SL bialg1-}
            \sigma\Big(\big(\mathrm{id}_{\frak g}\otimes\mathrm{id}_{\frak g}\otimes R_{\circ_{\frak g}}(x)\big){\bf SL}(r)\Big)=0,\;\forall x\in\frak g.
        \end{eqnarray}
        \item\label{1534} \eqref{eq:NS-coalg2} holds if and only if the following equation holds:
        \begin{eqnarray}\label{eq:SL bialg2-}
            \big(L_{\circ_{\frak g}}(x)\otimes\mathrm{id}_{\frak g}\otimes\mathrm{id}_{\frak g}
            +\mathrm{id}_{\frak g}\otimes L_{\circ_{\frak g}}(x)\otimes\mathrm{id}_{\frak g}
            +\mathrm{id}_{\frak g}\otimes\mathrm{id}_{\frak g}\otimes\mathrm{ad}_{\frak g}(x)\big){\bf SL}(r)=0,\;\forall x\in\frak g.
        \end{eqnarray}
        \item\label{1541} \eqref{eq:NS bialg1} holds automatically.
        \item\label{1542} \eqref{eq:NS bialg2} holds automatically.
    \end{enumerate}
\end{pro}
\begin{proof}
    \eqref{1530} It is clear that \eqref{eq:Lie coalg1} holds. Moreover, for all $x\in\frak g$, we have
    \begin{align*}
        \sigma(\mathrm{id}_{\frak g}\otimes\Delta)\Delta(x)
        &=\sigma\Big(\sum_i b_i\circ_{\frak g}x\otimes\Delta(a_i)-a_i\otimes\Delta(b_i\circ_{\frak g}x)\Big)\\
        &=\sigma\Big(\sum_{i,j}b_j\circ_{\frak g}a_i\otimes a_j\otimes b_i\circ_{\frak g}x-a_j\otimes b_j\circ_{\frak g}a_i\otimes b_i\circ_{\frak g}x\\
        &\quad-a_j\otimes a_i\otimes b_j\circ_{\frak g}(b_i\circ_{\frak g}x)+a_i\otimes a_j\otimes b_j\circ_{\frak g}(b_i\circ_{\frak g}x)\Big)\\
        &=\sigma\Big(\sum_{i,j}b_i\circ_{\frak g}a_j\otimes a_i\otimes b_j\circ_{\frak g}x-a_i\otimes b_i\circ_{\frak g}a_j\otimes b_j\circ_{\frak g}x\\
        &\quad-a_i\otimes a_j\otimes[b_i,b_j]_{\frak g}\circ_{\frak g}x\Big)\\
        &=-\sigma\Big(\big(\mathrm{id}_{\frak g}\otimes\mathrm{id}_{\frak g}\otimes R_{\circ_{\frak g}}(x)\big){\bf SL}(r)\Big).
    \end{align*}
    Hence \eqref{eq:Lie coalg2} holds if and only if \eqref{eq:SL bialg1-} holds.

    \eqref{1534} For all $x\in\frak g$, we have
    \begin{eqnarray*}
        &&(\Delta\otimes\mathrm{id}_{\frak g})\theta(x)-(\mathrm{id}_{\frak g}\otimes\theta)\theta(x)+(\tau\otimes\mathrm{id}_{\frak g})(\mathrm{id}_{\frak g}\otimes\theta)\theta(x)\\
        &&=\big(\Delta\otimes\mathrm{id}_{\frak g}-\mathrm{id}_{\frak g}\otimes\theta+(\tau\otimes\mathrm{id}_{\frak g})(\mathrm{id}_{\frak g}\otimes\theta)\big)
        \Big(\sum_i x\circ_{\frak g}a_i\otimes b_i+a_i\otimes[x,b_i]_{\frak g}\Big)\\
        &&=\sum_{i,j}b_j\circ_{\frak g}(x\circ_{\frak g}a_i)\otimes a_j\otimes b_i
        -a_j\otimes b_j\circ_{\frak g}(x\circ_{\frak g}a_i)\otimes b_i
        +b_j\circ_{\frak g}a_i\otimes a_j\otimes[x,b_i]_{\frak g}\\
        &&\ \ -a_j\otimes b_j\circ_{\frak g}a_i\otimes[x,b_i]_{\frak g}
        -x\circ_{\frak g}a_i\otimes b_i\circ_{\frak g}a_j\otimes b_j
        -x\circ_{\frak g}a_i\otimes a_j\otimes[b_i,b_j]_{\frak g}\\
        &&\ \ -a_i\otimes[x,b_i]_{\frak g}\circ_{\frak g}a_j\otimes b_j
        -a_i\otimes a_j\otimes[[x,b_i]_{\frak g},b_j]_{\frak g}
        +b_i\circ_{\frak g}a_j\otimes x\circ_{\frak g}a_i\otimes b_j\\
        &&\ \ +a_j\otimes x\circ_{\frak g}a_i\otimes[b_i,b_j]_{\frak g}
        +[x,b_i]_{\frak g}\circ_{\frak g}a_j\otimes a_i\otimes b_j
        +a_j\otimes a_i\otimes[[x,b_i]_{\frak g},b_j]_{\frak g}\\
        &&=\sum_{i,j}x\circ_{\frak g}(b_i\circ_{\frak g}a_j)\otimes a_i\otimes b_j
        -a_i\otimes x\circ_{\frak g}(b_i\circ_{\frak g}a_j)\otimes b_j
        +b_j\circ_{\frak g}a_i\otimes a_j\otimes[x,b_i]_{\frak g}\\
        &&\ \ -a_j\otimes b_j\circ_{\frak g}a_i\otimes[x,b_i]_{\frak g}
        -x\circ_{\frak g}a_i\otimes b_i\circ_{\frak g}a_j\otimes b_j
        -x\circ_{\frak g}a_i\otimes a_j\otimes[b_i,b_j]_{\frak g}\\
        &&\ \ +a_i\otimes a_j\otimes[x,[b_j,b_i]_{\frak g}]_{\frak g}
        +b_i\circ_{\frak g}a_j\otimes x\circ_{\frak g}a_i\otimes b_j
        +a_j\otimes x\circ_{\frak g}a_i\otimes[b_i,b_j]_{\frak g}\\
        &&=-\big(L_{\circ_{\frak g}}(x)\otimes\mathrm{id}_{\frak g}\otimes\mathrm{id}_{\frak g}
        +\mathrm{id}_{\frak g}\otimes L_{\circ_{\frak g}}(x)\otimes\mathrm{id}_{\frak g}
        +\mathrm{id}_{\frak g}\otimes\mathrm{id}_{\frak g}\otimes\mathrm{ad}_{\frak g}(x)\big){\bf SL}(r).
    \end{eqnarray*}
    Hence \eqref{eq:NS-coalg2} holds if and only if \eqref{eq:SL bialg2-} holds.

    \eqref{1541} For all $x,y\in\frak g$, we have
    \begin{align*}
        &\Delta(x\circ_{\frak g}y)-(\mathrm{id}_{\frak g}^{\otimes 2}-\tau)\Big(\big(L_{\circ_{\frak g}}(x)\otimes\mathrm{id}_{\frak g}\big)\Delta(y)+\big(\mathrm{id}_{\frak g}\otimes R_{\circ_{\frak g}}(y)\big)\theta(x)\Big)\\
        &=\sum_i b_i\circ_{\frak g}(x\circ_{\frak g}y)\otimes a_i-a_i\otimes b_i\circ_{\frak g}(x\circ_{\frak g}y)\\
        &\quad-(\mathrm{id}_{\frak g}^{\otimes 2}-\tau)\big(
        x\circ_{\frak g}(b_i\circ_{\frak g}y)\otimes a_i
        -x\circ_{\frak g}a_i\otimes b_i\circ_{\frak g}y
        +x\circ_{\frak g}a_i\otimes b_i\circ_{\frak g}y
        +a_i\otimes[x,b_i]_{\frak g}\circ_{\frak g}y\big)\\
        &=\sum_i b_i\circ_{\frak g}(x\circ_{\frak g}y)\otimes a_i-a_i\otimes b_i\circ_{\frak g}(x\circ_{\frak g}y)\\
        &\quad-x\circ_{\frak g}(b_i\circ_{\frak g}y)\otimes a_i
        +a_i\otimes x\circ_{\frak g}(b_i\circ_{\frak g}y)
        -a_i\otimes[x,b_i]_{\frak g}\circ_{\frak g}y
        +[x,b_i]_{\frak g}\circ_{\frak g}y\otimes a_i\\
        &=0.
    \end{align*}
    Hence \eqref{eq:NS bialg1} holds automatically.

    \eqref{1542} For all $x,y\in\frak g$, we have
    \begin{align*}
        &\theta([x,y]_{\frak g})
        -\big(L_{\circ_{\frak g}}(x)\otimes\mathrm{id}_{\frak g}+\mathrm{id}_{\frak g}\otimes\mathrm{ad}_{\frak g}(x)\big)\theta(y)
        +\big(L_{\circ_{\frak g}}(y)\otimes\mathrm{id}_{\frak g}+\mathrm{id}_{\frak g}\otimes\mathrm{ad}_{\frak g}(y)\big)\theta(x)\\
        &=\sum_i[x,y]_{\frak g}\circ_{\frak g}a_i\otimes b_i
        +a_i\otimes[[x,y]_{\frak g},b_i]_{\frak g}
        -x\circ_{\frak g}(y\circ_{\frak g}a_i)\otimes b_i
        -x\circ_{\frak g}a_i\otimes[y,b_i]_{\frak g}\\
        &\quad-y\circ_{\frak g}a_i\otimes[x,b_i]_{\frak g}
        -a_i\otimes[x,[y,b_i]_{\frak g}]_{\frak g}
        +y\circ_{\frak g}(x\circ_{\frak g}a_i)\otimes b_i
        +y\circ_{\frak g}a_i\otimes[x,b_i]_{\frak g}\\
        &\quad+x\circ_{\frak g}a_i\otimes[y,b_i]_{\frak g}
        +a_i\otimes[y,[x,b_i]_{\frak g}]_{\frak g}\\
        &=0.
    \end{align*}
    Hence \eqref{eq:NS bialg2} holds automatically.
\end{proof}

\begin{cor}\label{cor:4.3}
    Let $(\frak g,[-,-]_{\frak g},\circ_{\frak g})$ be a self-module Lie algebra and $r=\sum\limits_i a_i\otimes b_i\in\frak g\otimes\frak g$.
    Let $\Delta_r,\theta_r:\frak g\rightarrow\frak g\otimes\frak g$ be comultiplications given by \eqref{eq:NS com1} and \eqref{eq:NS com2}.
    Then $(\frak g,[-,-]_{\frak g},\circ_{\frak g},\Delta_r,\theta_r)$ is a self-module Lie bialgebra if and only if \eqref{eq:SL bialg1-} and \eqref{eq:SL bialg2-} hold. In particular,
    if $r$ is a solution of the SLYBE in $(\frak g,[-,-]_{\frak g},\circ_{\frak g})$, then $(\frak g,[-,-]_{\frak g},\circ_{\frak g},\Delta_r,\theta_r)$ is a self-module Lie bialgebra.
\end{cor}
\begin{proof}
    It follows directly from Proposition \ref{pro:4.2}.
\end{proof}

\begin{defi}\label{def:phi-condition}
    Let $\frak g$ be a vector space and let $\phi:\frak g\to\frak g$ be a linear isomorphism.
    We say that $r\in\frak g\otimes\frak g$ satisfies the {\bf $\phi$-condition} if
    \begin{equation}\label{eq:r}
        (\mathrm{id}_{\frak g}\otimes\phi)r+\tau(r)=0.
    \end{equation}
\end{defi}

Let $\frak g$ and $\frak h$ be vector spaces.
Then the linear bijection
$\frak h\otimes\frak g\cong\mathrm{Hom}(\frak h^*,\frak g)$
identifies any $r\in\frak h\otimes\frak g$ as a linear map
(still denoted by $r$) from $\frak h^*$ to $\frak g$ by
\begin{equation*}
    \langle r,\xi\otimes\eta\rangle
    =\langle r(\xi),\eta\rangle,
    \quad\forall\xi\in\frak h^*,\ \eta\in\frak g^*.
\end{equation*}

\begin{pro}\label{pro:2845}
    Let $(\frak g,[-,-]_{\frak g},\circ_{\frak g})$ be a self-module Lie algebra and $r\in\frak g\otimes\frak g$ be a solution of the SLYBE.
    Let $(\frak g,[-,-]_{\frak g},\circ_{\frak g},\Delta_r,\theta_r)$ be the self-module Lie bialgebra given by Corollary \ref{cor:4.3} and $\frak g\bowtie^{L_{\circ_{\frak g^*}}^*}_{L_{\circ_{\frak g}}^*}\frak g^*$ be the resulting Lie algebra structure on $\frak g\oplus\frak g^*$ in which the multiplications $[-,-]_{\frak g^*},\circ_{\frak g^*}:\frak g^*\otimes\frak g^*\rightarrow\frak g^*$ are given through \eqref{eq:2249}.
    Then $\Phi:\frak g\ltimes_{L_{\circ_{\frak g}}^*}\frak g^*\rightarrow\frak g\bowtie^{L_{\circ_{\frak g^*}}^*}_{L_{\circ_{\frak g}}^*}\frak g^*$ is an isomorphism of Lie algebras, where
    \begin{equation*}\label{eq:2851}
        \Phi(x)=x,\;\Phi(\xi )=-r(\xi )+\xi ,\;\forall x\in\frak g,\ \xi \in\frak g^*.
    \end{equation*}
    Moreover, $r$  satisfies the $\phi$-condition if and only if the following equation holds:
    \begin{equation}\label{eq:preserve}
        \mathcal B_{\phi}(x+\xi ,y+\eta )=\mathcal B_{\phi}\big(\Phi(x+\xi ),\Phi(y+\eta )\big),\;\forall x,y\in\frak g,\ \xi ,\eta \in\frak g^*.
    \end{equation}
\end{pro}
\begin{proof}
    Set $r=\sum\limits_i a_i\otimes b_i\in\frak g\otimes\frak g$.
    Let $x,y\in\frak g,\ \xi ,\eta \in\frak g^*$ and denote
    \begin{eqnarray*}
        \frak g\bowtie^{L_{\circ_{\frak g^*}}^*}_{L_{\circ_{\frak g}}^*}\frak g^*=(\frak g\oplus\frak g^*,[-,-]_{\bowtie}),\;
        \frak g\ltimes_{L_{\circ_{\frak g}}^*}\frak g^*=(\frak g\oplus\frak g^*,[-,-]_{\ltimes}).
    \end{eqnarray*}
    It is clear that $\Phi$ is a linear isomorphism and
    $[\Phi(x),\Phi(y)]_{\ltimes}=[\Phi(x),\Phi(y)]_{\bowtie}$.
    Moreover,
    \begin{align*}
        \Phi([x,\xi ]_{\ltimes})-[\Phi(x),\Phi(\xi )]_{\bowtie}
        &=\Phi(L_{\circ_{\frak g}}^*(x)\xi )-[x,-r(\xi )+\xi ]_{\bowtie}\\
        &=-r(L_{\circ_{\frak g}}^*(x)\xi )+L_{\circ_{\frak g}}^*(x)\xi +[x,r(\xi )]_{\frak g}-L_{\circ_{\frak g}}^*(x)\xi +L_{\circ_{\frak g^*}}^*(\xi )x\\
        &=-r(L_{\circ_{\frak g}}^*(x)\xi )+[x,r(\xi )]_{\frak g}+L_{\circ_{\frak g^*}}^*(\xi )x.
    \end{align*}
    Noting that
        \begin{align*}
            \langle-r(L_{\circ_{\frak g}}^*(x)\xi )+[x,r(\xi )]_{\frak g}+L_{\circ_{\frak g^*}}^*(\xi )x,\eta \rangle
            &=-\langle r(L_{\circ_{\frak g}}^*(x)\xi ),\eta \rangle-\langle r(\xi ),\mathrm{ad}_{\frak g}^*(x)\eta \rangle-\langle\theta_r(x),\xi \otimes \eta \rangle\\
            &=-\langle r,L_{\circ_{\frak g}}^*(x)\xi \otimes \eta \rangle-\langle r,\xi \otimes\mathrm{ad}_{\frak g}^*(x)\eta \rangle-\langle\theta_r(x),\xi \otimes \eta \rangle\\
            &=\langle \xi \otimes \eta ,(L_{\circ_{\frak g}}(x)\otimes\mathrm{id}_{\frak g}+\mathrm{id}_{\frak g}\otimes\mathrm{ad}_{\frak g}(x))r-\theta_r(x)\rangle\\
            &\overset{\eqref{eq:NS com2}}{=}0,
        \end{align*}
    we have $\Phi([x,\xi ]_{\ltimes})=[\Phi(x),\Phi(\xi )]_{\bowtie}$.
    Furthermore, since we have  ${\bf SL}(r)=0$, a straightforward computation shows that $ \Phi([\xi ,\eta ]_{\ltimes})=0=[\Phi(\xi ),\Phi(\eta )]_{\bowtie}$.
    Hence $\Phi$ is an isomorphism of Lie algebras.
    Furthermore, we have
    \begin{align*}
     &\mathcal B_{\phi}(x+\xi ,y+\eta )-\mathcal B_{\phi}\big(\Phi(x+\xi ),\Phi(y+\eta )\big) \\
        &=\mathcal B_{\phi}(x+\xi ,y+\eta )-\mathcal B_{\phi}\big(x-r(\xi )+\xi ,y-r(\eta )+\eta \big) =\mathcal B_{\phi}(r(\xi ),\eta )+\mathcal B_{\phi}(\xi ,r(\eta ))\\
        &=\langle\phi(r(\xi )),\eta \rangle+\langle r(\eta ),\xi \rangle =\langle r,\xi \otimes\phi^t(\eta )+\eta \otimes \xi \rangle=\langle \xi \otimes \eta ,(\mathrm{id}_{\frak
g}\otimes\phi)r+\tau(r)\rangle.
    \end{align*}
    Hence \eqref{eq:preserve} holds if and only if \eqref{eq:r} holds.
\end{proof}

\begin{defi}
    Let $(\frak g,[-,-]_{\frak g},\circ_{\frak g})$ be a self-module Lie algebra, $r\in\frak g\otimes\frak g$ be a solution of the SLYBE and $(\frak g,[-,-]_{\frak g},\circ_{\frak g},\Delta_r,\theta_r)$ be the self-module Lie bialgebra given in Corollary \ref{cor:4.3}.
    If $r$ satisfies the $\phi$-condition for a linear isomorphism $\phi$ on $\frak g$, then we say $(\frak g,[-,-]_{\frak g},\circ_{\frak g},\Delta_r,\theta_r)$ is a {\bf triangular self-module Lie bialgebra with respect to $\phi$}.
\end{defi}

Triangular Lie bialgebras, triangular pre-Lie bialgebras and
triangular anti-pre-Lie bialgebras are typical examples of
triangular self-module Lie bialgebras (Propositions \ref{2480}, \ref{2675} and \ref{2972}).

\subsection{$\mathcal{O}$-operators of self-module Lie algebras and self-bimodule pre-Lie algebras}

\begin{defi}
    Let $(\frak g,[-,-]_{\frak g},\circ_{\frak g})$ be a self-module Lie algebra with a representation $(\rho,l,r,V)$ and $T:V\rightarrow\frak g$ be a linear map.
    Suppose that $T$ is an $\mathcal{O}$-operator \cite{Kup} of $(\frak g,[-,-]_{\frak g})$ with respect to $(\rho,V)$, that is, the following equation holds:
    \begin{equation*}
        [Tu,Tv]_{\frak g}=T\big(\rho(Tu)v-\rho(Tv)u\big),\;\forall u,v\in V.\label{eq:SR-Lieop1}
    \end{equation*}
    If moreover the following equation holds:
    \begin{equation*}
        (Tu)\circ_{\frak g}(Tv)=T\big(l(Tu)v+r(Tv)u\big),\;\forall u,v\in V,\label{eq:SR-Lieop2}
    \end{equation*}
    then we say $T$ is an {\bf $\mathcal{O}$-operator of $(\frak g,[-,-]_{\frak g},\circ_{\frak g})$ with respect to $(\rho,l,r,V)$}.
\end{defi}

\begin{thm}\label{thm:3034}
    Let $(\frak g,[-,-]_{\frak g},\circ_{\frak g})$ be a self-module Lie algebra and $r\in\frak g\otimes\frak g$.
    \begin{enumerate}
        \item\label{3125} $r$ is a solution of the SLYBE in $(\frak g,[-,-]_{\frak g},\circ_{\frak g})$ if and only if $r$ is an $\mathcal{O}$-operator of $(\frak g,[-,-]_{\frak g})$ with respect to $(L_{\circ_{\frak g}}^*,\frak g^*)$, that is, the following equation holds:
        \begin{eqnarray}\label{eq:Lie op ybe}
            [r(\xi ),r(\eta )]_{\frak g}=r\Big(L_{\circ_{\frak g}}^*\big(r(\xi )\big)\eta -L_{\circ_{\frak g}}^*\big(r(\eta )\big)\xi \Big),\;\forall \xi ,\eta \in\frak g^*.
        \end{eqnarray}
        \item\label{1655} Suppose that $r$ satisfies the $\phi$-condition for a linear isomorphism $\phi$ on $\frak g$. Let $(\frak g,[-,-]_{\frak g},\circ_\phi)$ be the $\phi$-conjugate self-module Lie algebra of $(\frak g,[-,-]_{\frak g},\circ_{\frak g})$. Then $r$ is a solution of the SLYBE in $(\frak g,[-,-]_{\frak g},\circ_{\frak g})$ if and only if the following equation holds:
        \begin{eqnarray}\label{eq:SR-Lie op ybe1}
            r(\xi )\circ_\phi r(\eta )=r\Big(\mathrm{ad}_{\frak g}^*\big(r(\xi )\big)\eta +R_{\circ_{\frak g}}^*[1,\phi,1]\big(r(\eta )\big)\xi \Big),\;\forall \xi ,\eta \in\frak g^*.
        \end{eqnarray}
        In this case, $r$ is an $\mathcal{O}$-operator of the self-module Lie algebra $(\frak g,[-,-]_{\frak g},\circ_\phi)$ with respect to $ (L_{\circ_{\frak g}}^*,\mathrm{ad}_{\frak g}^*,R_{\circ_{\frak g}}^*[1,\phi,1],\frak g^* )$.
    \end{enumerate}
\end{thm}
\begin{proof}
    Set $r=\sum\limits_i a_i\otimes b_i\in\frak g\otimes\frak g$.
    Let $\xi ,\eta ,\gamma \in\frak g^*$. Then we have
    \begin{align}
        \langle[r(\xi ),r(\eta )]_{\frak g},\gamma \rangle
        &=\langle r(\xi ),\mathrm{ad}_{\frak g}^*\big(r(\eta )\big)\gamma \rangle
        =\sum_i\langle a_i,\xi \rangle\langle b_i,\mathrm{ad}_{\frak g}^*\big(r(\eta )\big)\gamma \rangle\nonumber\\
        &=-\sum_i\langle a_i,\xi \rangle\langle r(\eta ),\mathrm{ad}_{\frak g}^*(b_i)\gamma \rangle
        =-\sum_{i,j}\langle a_i,\xi \rangle\langle a_j,\eta \rangle\langle b_j,\mathrm{ad}_{\frak g}^*(b_i)\gamma \rangle\nonumber\\
        &=\sum_{i,j}\langle a_i\otimes a_j\otimes[b_i,b_j]_{\frak g},\xi \otimes \eta \otimes \gamma \rangle,\label{eq:76}
    \end{align}
    and similarly
    \begin{align*}
        -\langle r\Big(L_{\circ_{\frak g}}^*\big(r(\xi )\big)\eta \Big),\gamma \rangle
        &=\sum_{i,j}\langle a_j\otimes b_j\circ_{\frak g}a_i\otimes b_i,\xi \otimes \eta \otimes \gamma \rangle,\\
        \langle r\Big(L_{\circ_{\frak g}}^*\big(r(\eta )\big)\xi \Big),\gamma \rangle
        &=-\sum_{i,j}\langle b_j\circ_{\frak g}a_i\otimes a_j\otimes b_i,\xi \otimes \eta \otimes \gamma \rangle.
    \end{align*}
    Thus ${\bf SL}(r)=0$ in $(\frak g,[-,-]_{\frak g},\circ_{\frak g})$ if and only if \eqref{eq:Lie op ybe} holds. Therefore we get Item \eqref{3125}.
    Moreover, replacing the multiplication $[-,-]_{\frak g}$ by $\circ_\phi$ in \eqref{eq:76}, we have
    \begin{align*}
        \langle r(\xi )\circ_\phi r(\eta ),\gamma \rangle
        &= \sum_{i,j}\langle a_i\otimes a_j\otimes b_i\circ_\phi b_j,\xi \otimes \eta \otimes \gamma \rangle\\
        &= \sum_{i,j}\langle(\phi^{-1}\otimes\mathrm{id}_{\frak g}\otimes\mathrm{id}_{\frak g})\big(b_i\circ_{\frak g}\phi(b_j)\otimes a_i\otimes a_j\big),\gamma \otimes \xi \otimes \eta \rangle\\
        &\overset{\eqref{eq:r}}{=} -\sum_{i,j}\langle(\phi^{-1}\otimes\mathrm{id}_{\frak g}\otimes\mathrm{id}_{\frak g})\big(b_i\circ_{\frak g}a_j\otimes a_i\otimes b_j\big),\gamma \otimes \xi \otimes \eta \rangle.
    \end{align*}
    Similarly, we have
    \begin{align*}
        \langle r\Big(\mathrm{ad}_{\frak g}^*\big(r(\xi )\big)\eta \Big),\gamma \rangle
        &=-\sum_{i,j}\langle(\phi^{-1}\otimes\mathrm{id}_{\frak g}\otimes\mathrm{id}_{\frak g})\big(a_i\otimes a_j\otimes[b_i,b_j]_{\frak g}\big),\gamma \otimes \xi \otimes \eta \rangle,\\
        \langle r\Big(R_{\circ_{\frak g}}^*[1,\phi,1]\big(r(\eta )\big)\xi \Big),\gamma \rangle
        &=-\sum_{i,j}\langle(\phi^{-1}\otimes\mathrm{id}_{\frak g}\otimes\mathrm{id}_{\frak g})\big(a_i\otimes b_i\circ_{\frak g}a_j\otimes b_j\big),\gamma \otimes \xi \otimes \eta \rangle.
    \end{align*}
    Hence we have
    \begin{equation*}
    \langle r(\xi )\circ_\phi r(\eta )-r\Big(\mathrm{ad}_{\frak g}^*\big(r(\xi )\big)\eta +R_{\circ_{\frak g}}^*[1,\phi,1]\big(r(\eta )\big)\xi \Big),\gamma \rangle =\langle(\phi^{-1}\otimes\mathrm{id}_{\frak g}\otimes\mathrm{id}_{\frak g}){\bf SL}(r),\gamma \otimes \xi \otimes \eta \rangle.
    \end{equation*}
    Thus ${\bf SL}(r)=0$ in $(\frak g,[-,-]_{\frak g},\circ_{\frak g})$ if and only if \eqref{eq:SR-Lie op ybe1} holds.
    Hence we get Item \eqref{1655}.
\end{proof}

\begin{pro}\label{pro:3094}
    Let $(\frak g,[-,-]_{\frak g},\circ_{\frak g})$ be a self-module Lie algebra and $r\in\frak g\otimes\frak g$ be nondegenerate which is identified as an invertible linear map $r:\frak g^*\rightarrow\frak g$.
    \begin{enumerate}
        \item\label{1789} Let $\phi:\frak g\rightarrow\frak g$ be the linear isomorphism given by
        \begin{equation}\label{1791}
            \phi(x)=-r^t r^{-1}(x),\;\forall x\in\frak g.
        \end{equation}
        Then $r$ satisfies the $\phi$-condition.
        \item\label{1794} Let $\omega$ be the nondegenerate bilinear form on $\frak g$ given by
        \begin{equation}\label{1901}
            \omega(x,y)=\langle r^{-1}(x),y\rangle,\;\forall x,y\in\frak g.
        \end{equation}
        Then $r:\frak g^*\rightarrow\frak g$ is an $\mathcal{O}$-operator of $(\frak g,[-,-]_{\frak g})$ with respect to $(L_{\circ_{\frak g}}^*,\frak g^*)$ if and only if the following equation holds:
        \begin{eqnarray}\label{eq:2-coc}
            \omega([x,y]_{\frak g},z)=\omega\big(x,y\circ_{\frak g}z\big)-\omega\big(y,x\circ_{\frak g}z\big),\;\forall x,y,z\in\frak g.
        \end{eqnarray}
    \end{enumerate}
\end{pro}
\begin{proof}
    For all $\xi ,\eta \in\frak g^*$, we have
    \begin{equation*}
        \langle(\mathrm{id}_{\frak g}\otimes\phi)r,\xi \otimes \eta \rangle
         =\langle\phi r(\xi ),\eta \rangle
        \overset{\eqref{1791}}{=}-\langle r^t(\xi ),\eta \rangle
     =-\langle r(\eta ),\xi \rangle
        =-\langle\tau(r),\xi \otimes \eta \rangle.
    \end{equation*}
    Hence we get Item \eqref{1789}.
    Moreover, let $x,y,z\in\frak g$ and $\xi =\omega^\natural(x),\ \eta =\omega^\natural(y)$. By \eqref{1901}, we have $r(\xi )=x,\ r(\eta )=y$. Furthermore, we have
    \begin{align*}
        \langle r^{-1}[r(\xi ),r(\eta )]_{\frak g},z\rangle
        &=\langle\omega^\natural([x,y]_{\frak g}),z\rangle=\omega([x,y]_{\frak g},z),\\
        \langle L_{\circ_{\frak g}}^*\big(r(\xi )\big)\eta ,z\rangle
        &=-\langle \eta ,r(\xi )\circ_{\frak g}z\rangle
        =-\langle\omega^\natural(y),x\circ_{\frak g}z\rangle
        =-\omega(y,x\circ_{\frak g}z),\\
        -\langle L_{\circ_{\frak g}}^*\big(r(\eta )\big)\xi ,z\rangle
        &=\omega(x,y\circ_{\frak g}z).
    \end{align*}
    Hence \eqref{eq:2-coc} holds if and only if \eqref{eq:Lie op ybe} holds. Hence we get Item \eqref{1794}.
\end{proof}

We next study the solutions of the SLYBE in semi-direct product self-module Lie algebras.

\begin{pro}\label{pro:4.7}
    Let $\frak g=(\frak g,[-,-]_{\frak g},\circ_{\frak g})$ be a self-module Lie algebra, $\phi:\frak g\rightarrow\frak g$ be a linear isomorphism and $\frak g_\phi=(\frak g,[-,-]_{\frak g},\circ_\phi)$ be the $\phi$-conjugate self-module Lie algebra.
    Let $(\rho,l,r,V)$ be a representation of $(\frak g,[-,-]_{\frak g},\circ_{\frak g})$, $f:V\rightarrow V$ be a linear isomorphism and
    $$\big(l^*,\rho^*[(f^{-1})^t,1,f^t],r^*[(f^{-1})^t,\phi,1],V^*\big)$$
    be the corresponding representation of $\frak g_\phi$ in Proposition~\ref{848}.
    Let
    $$(\widehat{\frak g},[-,-]_{\widehat{\frak g}},\circ_{\widehat{\frak g}})
    =\frak g_\phi\ltimes_{l^*,\rho^*[(f^{-1})^t,1,f^t],r^*[(f^{-1})^t,\phi,1]}V^*$$
    be the resulting semi-direct product self-module Lie algebra structure on $\widehat{\frak g}=\frak g\oplus V^*$.
    Let $T,S:V\rightarrow\frak g$ be linear maps which are identified as elements in $V^*\otimes\frak g\subset\widehat{\frak g}\otimes\widehat{\frak g}$.
    Then $r=T+\tau(S)$ is a solution of the SLYBE in $(\widehat{\frak g},[-,-]_{\widehat{\frak g}},\circ_{\widehat{\frak g}})$ if and only if the following conditions hold:
    \begin{align*}
        [Tu,Tv]_{\frak g}&=T\big(f\rho(Tu)f^{-1}(v)-f\rho(Tv)f^{-1}(u)\big),\\
        (Tu)\circ_\phi(Sv)&=S\Big(l(Tu)v-r\phi(Sv)f^{-1}(u)\Big),\;\forall u,v\in V.
    \end{align*}
    Moreover, suppose that $P_{\frak g}:\frak g\rightarrow\frak g$ and $P_V:V\rightarrow V$ are linear isomorphisms. Then $r$ satisfies the $(P_{\frak g}+P_V^t)$-condition if and only if the following equations hold:
    \begin{eqnarray}
        P_{\frak g}T+S=0,&&\label{1838}\\
        T+SP_V=0.&&\label{1839}
    \end{eqnarray}
\end{pro}
\begin{proof}
    The first assertion follows by a tensor computation similar to that in \cite[Section 2]{Bai2007}.
    Let $\{v_1,\cdots,v_n\}$ be a basis of $V$ and $\{v_1^*,\cdots,v_n^*\}$ be the dual basis.
    Then we have
    \begin{align*}
        &\big(\mathrm{id}_{\widehat{\frak g}}\otimes(P_{\frak g}+P_V^t)\big)\big(T+\tau(S)\big)+\tau(T)+S\\
        &=\sum_i v_i^*\otimes P_{\frak g}(Tv_i)+S(v_i)\otimes P_V^t(v_i^*)+T(v_i)\otimes v_i^*+v_i^*\otimes S(v_i)\\
        &=\sum_i v_i^*\otimes(P_{\frak g}T+S)v_i+(SP_V+T)v_i\otimes v_i^*.
    \end{align*}
    Thus $r$ satisfies the $(P_{\frak g}+P_V^t)$-condition if and only if \eqref{1838} and \eqref{1839} hold.
\end{proof}

We next consider special cases for which the results in Proposition \ref{pro:4.7} can be made simpler.

\begin{cor}\label{cor:4.16-}
    Let $(\frak g,[-,-]_{\frak g},\circ_{\frak g})$ be a self-module Lie algebra, $(\rho,l,r,V)$ be a representation of $(\frak g,[-,-]_{\frak g},\circ_{\frak g})$ and
    \begin{equation*}
        (\widehat{\frak g},[-,-]_{\widehat{\frak g}},\circ_{\widehat{\frak g}})
        =\frak g\ltimes_{l^*,\rho^*,-r^*}V^*,\;
        (\widehat{\frak g},[-,-]_{\widehat{\frak g}},\bullet_{\widehat{\frak g}})
        =\frak g\ltimes_{l^*,\rho^*,r^*}V^*
    \end{equation*}
    be the corresponding semi-direct product self-module Lie algebras.
    Let $T:V\rightarrow\frak g$ be a linear map which is identified as an element in $V^*\otimes\frak g\subset\widehat{\frak g}\otimes\widehat{\frak g}$.
    Then the following statements are equivalent:
    \begin{enumerate}
        \item\label{4302} $T$ is an $\mathcal{O}$-operator of $(\frak g,[-,-]_{\frak g},\circ_{\frak g})$ with respect to $(\rho,l,r,V)$.
        \item\label{4304} $s=T+\tau(T)$ is a symmetric solution of the SLYBE in $(\widehat{\frak g},[-,-]_{\widehat{\frak g}},\circ_{\widehat{\frak g}})$.
        \item\label{4306} $r=T-\tau(T)$ is an antisymmetric solution of the SLYBE in $(\widehat{\frak g},[-,-]_{\widehat{\frak g}},\bullet_{\widehat{\frak g}})$.
    \end{enumerate}
\end{cor}
\begin{proof}
    The equivalence \eqref{4302}$\Longleftrightarrow$\eqref{4304} is a direct consequence of Proposition \ref{pro:4.7} by taking $\phi=\mathrm{id}_{\frak g},f=-\mathrm{id}_V,S=T$, and the equivalence \eqref{4302}$\Longleftrightarrow$\eqref{4306} is a direct consequence of Proposition \ref{pro:4.7} by taking $\phi=\mathrm{id}_{\frak g},f=\mathrm{id}_V,S=-T$.
\end{proof}

\begin{rmk}
    There are other ways to simplify Proposition \ref{pro:4.7}. For example, taking $T=0$ in Proposition \ref{pro:4.7}, we know that $\tau(S)$ is a solution of the SLYBE in $(\widehat{\frak g},[-,-]_{\widehat{\frak g}},\circ_{\widehat{\frak g}})$ if and only if the following equation holds:
    \begin{eqnarray}\label{3863}
        S\Big(r\big(\phi(Sv)\big)u\Big)=0,\;\forall u,v\in V.
    \end{eqnarray}
    In particular, if the right action $r$ is zero, then \eqref{3863}
    holds for every linear map $S:V\to\frak g$. Thus $\tau(S)$ is a
    solution of the SLYBE and gives rise to a self-module Lie bialgebra.
    For $S\ne 0$, note that the tensor $\tau(S)$ does not satisfy the $P$-condition
    for any linear isomorphism $P$ on $\widehat{\frak g}$, since the two terms
    in $(\mathrm{id}_{\widehat{\frak g}}\otimes P)\tau(S)+S$ lie in
    $\frak g\otimes\widehat{\frak g}$ and
    $V^*\otimes\widehat{\frak g}$ respectively, and cannot cancel.
\end{rmk}

Recall that a {\bf pre-Lie algebra} \cite{Bai2021.2,Bur} is a vector space $\frak g$ with a multiplication $\cdot_{\frak g}:\frak g\otimes\frak g\rightarrow\frak g$ such that the following equation holds:
\begin{equation*}\label{eq:pre-Lie}
    (x\cdot_{\frak g}y)\cdot_{\frak g}z-x\cdot_{\frak g}(y\cdot_{\frak g}z)
    =(y\cdot_{\frak g}x)\cdot_{\frak g}z-y\cdot_{\frak g}(x\cdot_{\frak g}z),\;\forall x,y,z\in\frak g.
\end{equation*}
For a pre-Lie algebra $(\frak g,\cdot_{\frak g})$, there is a Lie algebra $(\frak g,[-,-]_{\frak g})$ given by
\begin{equation}\label{eq:sub-adj}
    [x,y]_{\frak g}=x\cdot_{\frak g}y-y\cdot_{\frak g}x,\;\forall x,y\in\frak g.
\end{equation}
In this case, we say that $(\frak g,[-,-]_{\frak g})$ is the {\bf sub-adjacent Lie algebra} of $(\frak g,\cdot_{\frak g})$ and $(\frak g,\cdot_{\frak g})$ is a {\bf compatible pre-Lie algebra} of $(\frak g,[-,-]_{\frak g})$.
Moreover, $(L_{\cdot_{\frak g}},\frak g)$ is a representation of $(\frak g,[-,-]_{\frak g})$.
Now we introduce the notion of a self-bimodule pre-Lie algebra.

\begin{defi}
    Let $(\frak g,\cdot_{\frak g})$ be a pre-Lie algebra. If there are multiplications $\succ_{\frak g},\prec_{\frak g}:\frak g\otimes\frak g\rightarrow\frak g$ satisfying the following equations:
    \begin{align*}
        (x\cdot_{\frak g}y)\succ_{\frak g}z-x\succ_{\frak g}(y\succ_{\frak g}z)
        &=(y\cdot_{\frak g}x)\succ_{\frak g}z-y\succ_{\frak g}(x\succ_{\frak g}z),\\
        x\succ_{\frak g}(y\prec_{\frak g}z)-(x\succ_{\frak g}y)\prec_{\frak g}z
        &=y\prec_{\frak g}(x\cdot_{\frak g}z)-(y\prec_{\frak g}x)\prec_{\frak g}z,\;\forall x,y,z\in\frak g,
    \end{align*}
    then we say $(\frak g,\cdot_{\frak g},\succ_{\frak g},\prec_{\frak g})$ is a {\bf self-bimodule pre-Lie algebra}.
\end{defi}

Recall \cite{Bai2008} that a {\bf representation} of a pre-Lie
algebra $(\frak g,\cdot_{\frak g})$ is a triple $(l,r,V)$, where
$V$ is a vector space and $l,r:\frak
g\rightarrow\mathrm{End}_{\mathbb K}(V)$ are linear maps
satisfying the following equations:
\begin{align*}
    l(x\cdot_{\frak g}y)v-l(x)l(y)v&=l(y\cdot_{\frak g}x)v-l(y)l(x)v,\\
    r(y)l(x)v-l(x)r(y)v&=r(y)r(x)v-r(x\cdot_{\frak g}y)v,\;\forall x,y\in\frak g,\ v\in V.
\end{align*}
In fact, $(l,r,V)$ is a representation of $(\frak g,\cdot_{\frak g})$ if and only if there is a pre-Lie algebra structure $(\frak g\oplus V,\cdot_d)=\frak g\ltimes_{l,r}V$ on $\frak g\oplus V$ with $\cdot_d$ given by
\begin{eqnarray*}
    (x+u)\cdot_d(y+v)=x\cdot_{\frak g}y+l(x)v+r(y)u,\;\forall x,y\in\frak g,\ u,v\in V.
\end{eqnarray*}
The following characterization justifies the terminology of a self-bimodule pre-Lie algebra.

\begin{pro}\label{4894}
    Let $(\frak g,\cdot_{\frak g})$ be a pre-Lie algebra. Then there exist multiplications $\succ_{\frak g},\prec_{\frak g}:\frak g\otimes\frak g\rightarrow\frak g$ such that $(\frak g,\cdot_{\frak g},\succ_{\frak g},\prec_{\frak g})$ is a self-bimodule pre-Lie algebra if and only if there is a representation $(l,r,\frak g)$ of $(\frak g,\cdot_{\frak g})$. Here the mutual relation is given by $l=L_{\succ_{\frak g}}$ and $r=R_{\prec_{\frak g}}$, or equivalently,
    \begin{equation*}
        x\succ_{\frak g}y=l(x)y,\;
        x\prec_{\frak g}y=r(y)x,\;\forall x,y\in\frak g.
    \end{equation*}
\end{pro}
\begin{proof}
    It follows by a straightforward verification.
\end{proof}

\begin{rmk}
As in Remark \ref{rmk530}, one can check that there is a one-to-one correspondence between self-bimodule pre-Lie algebras and NS-pre-Lie algebras \cite{GZ2023}.
\end{rmk}

\begin{pro}\label{pro:pre-Ns-}
    Let $(\frak g,[-,-]_{\frak g},\circ_{\frak g})$ be a self-module Lie algebra with a representation $(\rho,l,r,V)$.
    Let $T:V\rightarrow\frak g$ be an $\mathcal{O}$-operator of $(\frak g,[-,-]_{\frak g},\circ_{\frak g})$ with respect to $(\rho,l,r,V)$.
    Then there is a self-bimodule pre-Lie algebra structure $(V,\cdot_V,\succ_V,\prec_V)$ on $V$, where the multiplications $\cdot_V,\succ_V$ and $\prec_V$ are given by
    \begin{equation*}
        u\cdot_V v=\rho(Tu)v,\;
        u\succ_V v=l(Tu)v,\;
        u\prec_V v=-r(Tu)v,\;\forall u,v\in V.
    \end{equation*}
    Conversely, let $(\frak g,\cdot_{\frak g},\succ_{\frak g},\prec_{\frak g})$ be a self-bimodule pre-Lie algebra.
    Then there is a self-module Lie algebra $(\frak g,[-,-]_{\frak g},\circ_{\frak g})$ with $[-,-]_{\frak g}$ given by \eqref{eq:sub-adj} and $\circ_{\frak g}$ given by
    \begin{equation*}
        x\circ_{\frak g}y=x\succ_{\frak g}y-y\prec_{\frak g}x,\;\forall x,y\in\frak g,\label{eq:3178}
    \end{equation*}
    which is called the {\bf sub-adjacent self-module Lie algebra} of $(\frak g,\cdot_{\frak g},\succ_{\frak g},\prec_{\frak g})$.
    Moreover, $(L_{\cdot_{\frak g}},L_{\succ_{\frak g}}$,
    $-L_{\prec_{\frak g}},\frak g)$ is a representation of $(\frak g,[-,-]_{\frak g},\circ_{\frak g})$, and the identity map $\mathrm{id}_{\frak g}$ is an $\mathcal{O}$-operator of $(\frak g$,
    $[-,-]_{\frak g},\circ_{\frak g})$ with respect to $(L_{\cdot_{\frak g}},L_{\succ_{\frak g}},-L_{\prec_{\frak g}},\frak g)$.
\end{pro}
\begin{proof}
    It follows from a straightforward verification.
\end{proof}

\begin{cor}\label{2110}
    Let $(\frak g,\cdot_{\frak g},\succ_{\frak g},\prec_{\frak g})$ be a self-bimodule pre-Lie algebra and $(\frak g,[-,-]_{\frak g},\circ_{\frak g})$ be the sub-adjacent self-module Lie algebra.
    Let $\{e_1,\cdots,e_n\}$ be a basis of $\frak g$ and $\{e_1^*,\cdots,e_n^*\}$ be the dual basis.
    \begin{enumerate}
        \item
        $s=\sum\limits_{i=1}^n e_i^*\otimes e_i+e_i\otimes e_i^*$
        is a symmetric solution of the SLYBE in the self-module Lie algebra
        $\frak g\ltimes_{L_{\succ_{\frak g}}^*,L_{\cdot_{\frak g}}^*,L_{\prec_{\frak g}}^*}\frak g^*$.
        \item
        $r=\sum\limits_{i=1}^n e_i^*\otimes e_i-e_i\otimes e_i^*$
        is an antisymmetric solution of the SLYBE in the self-module Lie algebra
        $\frak g\ltimes_{L_{\succ_{\frak g}}^*,L_{\cdot_{\frak g}}^*,-L_{\prec_{\frak g}}^*}\frak g^*$.
    \end{enumerate}
\end{cor}
\begin{proof}
    The conclusions follow from Corollary \ref{cor:4.16-} and Proposition \ref{pro:pre-Ns-} by taking $(\rho,l,r,V)$
    $=
    \big(L_{\cdot_{\frak g}},L_{\succ_{\frak g}},-L_{\prec_{\frak g}},\frak g\big)$ and observing that $\mathrm{id}_{\frak g}=\sum\limits_{i=1}^n e_i^*\otimes e_i$.
\end{proof}

The following Figure \ref{fig:main-correspondences} summarizes the main correspondences and constructions in Sections \ref{sec:3} and \ref{sec:4} under the hypotheses of the cited results.

\begin{figure}[htbp]
    \centering
    \[
    \xymatrix@C=2.5cm@R=1.3cm{
        \txt{Standard generic Manin triples\\of Lie algebras}
        \ar@{<=>}[r]^-{\mathrm{Lem.}\ \ref{1452}}
        \ar@{<=>}[d]_-{\mathrm{Prop.}\ \ref{1559}}
        &
        \txt{Matched pairs of Lie algebras}
        \ar@{<=>}[d]^-{\mathrm{Prop.}\ \ref{pro:equ2}}
        \\
        \txt{Families of matched pairs of\\
            conjugate self-module Lie algebras}
        \ar@{<=>}[r]^-{\mathrm{Thm.}\ \ref{cor:2532}}
        &
        \txt{Self-module Lie bialgebras}
        \\
        \txt{$\mathcal O$-operators of Lie algebras}
        \ar@{<=>}[r]^-{\mathrm{Thm.}\ \ref{thm:3034}}
        &
        \txt{Solutions of the SLYBE}
        \ar[u]_-{\mathrm{Cor.}\ \ref{cor:4.3}}
        \\
        \txt{Self-bimodule pre-Lie algebras}
        \ar[r]^-{\mathrm{Prop.}\ \ref{pro:pre-Ns-}}
        \ar[ur]^-{\mathrm{Cor.}\ \ref{2110}}
        &
        \txt{$\mathcal O$-operators of\\self-module Lie algebras}
        \ar[u]_-{\mathrm{Cor.}\ \ref{cor:4.16-}}
    }
    \]
    \caption{Main correspondences and constructions for self-module Lie bialgebras.}
    \label{fig:main-correspondences}
\end{figure}
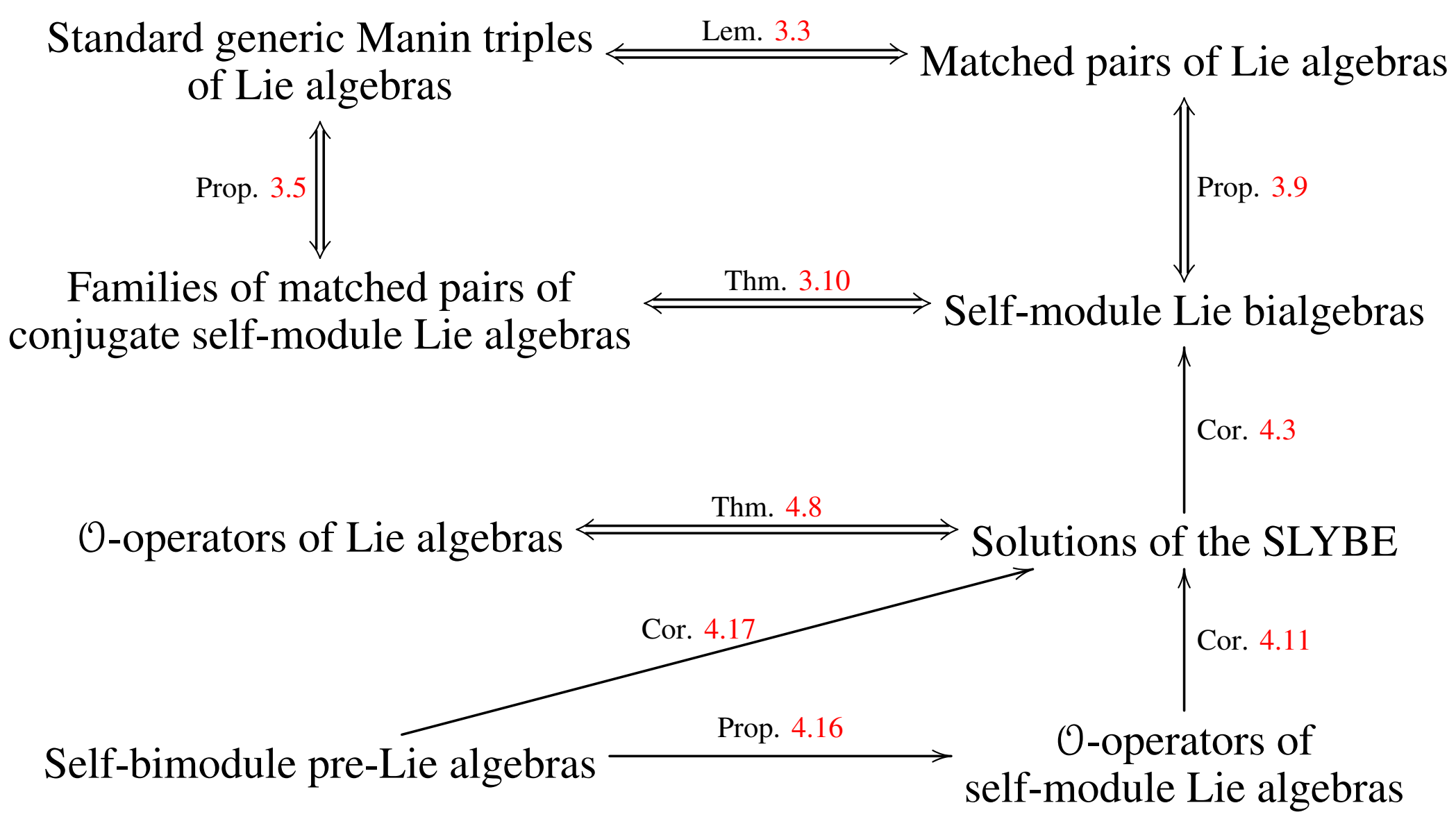

\begin{ex}\label{ex:nonclassical}
    Continue with Example \ref{557}.
    Then $(\frak n,\circ_{\frak n})$ is a pre-Lie algebra
    whose sub-adjacent Lie algebra is $(\frak n, [-,-]_{\frak n})$,
    and $(\frak n,\circ_{\frak n},\circ_{\frak n},\circ_{\frak n})$ is a self-bimodule pre-Lie algebra.
    Let $\{e_1^*,e_2^*\}$ be the dual basis and set $\frak h=\frak n\oplus\frak n^*$.
    By Corollary \ref{2110},
    \begin{equation*}
        s=e_1^*\otimes e_1+e_1\otimes e_1^*
        +e_2^*\otimes e_2+e_2\otimes e_2^*
    \end{equation*}
    is a symmetric solution of the SLYBE in the self-module Lie algebra
    \begin{equation*}
        (\frak h,[-,-]_{\frak h},\circ_1)
        =\frak n\ltimes_{L_{\circ_{\frak n}}^*,L_{\circ_{\frak n}}^*,L_{\circ_{\frak n}}^*}\frak n^*,
    \end{equation*}
    where the nonzero Lie brackets on $(\frak h,[-,-]_{\frak h})$ are given by \eqref{572} and
    \begin{eqnarray*}
        [e_1,e_1^*]_{\frak h}=-e_2^*,\;
        [e_2,e_2^*]_{\frak h}=-e_2^*,
    \end{eqnarray*}
    and the nonzero products of $\circ_1$ are given by
    \begin{eqnarray*}
        &&e_1\circ_1e_2=e_1,\;
        e_2\circ_1e_1=-e_1,\\
        &&e_1\circ_1e_1^*=e_1^*\circ_1e_1=-e_2^*,\;
        e_2\circ_1e_2^*=e_2^*\circ_1e_2=-e_2^*.
    \end{eqnarray*}
    By Corollary \ref{cor:4.3}, there is a triangular self-module Lie bialgebra $(\frak h,[-,-]_{\frak h},\circ_1,\Delta_s,\theta_s)$ with respect to $-\mathrm{id}_{\frak h}$, where the nonzero coproducts are given by
    \begin{eqnarray*}
        &&\Delta_s(e_2)=e_1\otimes e_1^*-e_1^*\otimes e_1
        +e_2\otimes e_2^*-e_2^*\otimes e_2,\\
        &&\theta_s(e_2)=-s,\;
        \Delta_s(e_1^*)=\theta_s(e_1^*)
        =e_1^*\otimes e_2^*-e_2^*\otimes e_1^*.
    \end{eqnarray*}
 It is straightforward to show that this self-module Lie bialgebra does not give a Lie bialgebra, or a pre-Lie bialgebra, or an anti-pre-Lie bialgebra (see
 Remark~\ref{rmk:c-ex}).
    Similarly, by Corollary \ref{2110} again,
    \begin{eqnarray*}
        r=e_1^*\otimes e_1-e_1\otimes e_1^*
        +e_2^*\otimes e_2-e_2\otimes e_2^*
    \end{eqnarray*}
    is an antisymmetric solution of the SLYBE in the self-module Lie algebra
    \begin{eqnarray*}
        (\frak h,[-,-]_{\frak h},\circ_2)
        =\frak n\ltimes_{L_{\circ_{\frak n}}^*,L_{\circ_{\frak n}}^*,-L_{\circ_{\frak n}}^*}\frak n^*,
    \end{eqnarray*}
    where $\circ_2=[-,-]_{\frak h}$.
    Consequently, by Corollary \ref{cor:4.3}, there is a triangular self-module Lie bialgebra $(\frak h,[-,-]_{\frak h},\circ_2,\Delta_r,\theta_r)$ with respect to $\mathrm{id}_{\frak h}$, where the nonzero coproducts are given by
    \begin{eqnarray*}
        &&\Delta_r(e_2)=\theta_r(e_2)=-r,\\
        &&\Delta_r(e_1^*)=\theta_r(e_1^*)
        =e_1^*\otimes e_2^*-e_2^*\otimes e_1^*.
    \end{eqnarray*}
   It is straightforward to show that this self-module Lie bialgebra gives rise to a triangular  Lie bialgebra (see
 Remark~\ref{rmk:c-ex}).
\end{ex}

\section{Classical examples of self-module Lie bialgebras}\label{2584}\
We show that Lie bialgebras, pre-Lie bialgebras and anti-pre-Lie
bialgebras are examples of self-module Lie bialgebras. See
Table~\ref{tab:classical-specializations} for an overview. Most statements follow by direct verification, so their proofs are omitted.

\subsection{Lie bialgebras as a classical example of self-module Lie bialgebras}\

\begin{ex}\label{3234}
Let $(\frak g,[-,-]_{\frak g})$ be a Lie algebra. Then $(\mathrm{ad}_{\frak g},\frak g)$ is a representation of $(\frak g,[-,-]_{\frak g})$.
By Proposition \ref{473},  $(\frak g,[-,-]_{\frak g},[-,-]_{\frak g})$ is a self-module Lie algebra.
\end{ex}

Recall that a {\bf quadratic Lie algebra} $(\frak g,[-,-]_{\frak g},\mathcal{B})$ is a Lie algebra $(\frak g,[-,-]_{\frak g} )$ together with a nondegenerate symmetric bilinear form $\mathcal{B}$ which is {\bf invariant} in the sense that
\begin{equation*}
\mathcal{B}([x,y]_{\frak g},z)=\mathcal{B}(x,[y,z]_{\frak g}),\;\forall x,y,z\in \frak g.
\end{equation*}

\begin{pro}\label{2887}
Let $(\frak g,[-,-]_{\frak g},\mathcal{B})$ be a quadratic Lie algebra.
Then the self-module Lie algebra induced from $(\frak g,[-,-]_{\frak g},\mathcal{B})$ is exactly $(\frak g,[-,-]_{\frak g},[-,-]_{\frak g})$ in Example \ref{3234}.
\end{pro}

\begin{defi}\label{2388}
A {\bf standard Manin triple of Lie algebras with respect to  the nondegenerate symmetric invariant bilinear form} \cite{Cha} is a standard generic Manin triple of Lie algebras $\big((\frak g\oplus\frak g^*,[-,-]_d,\mathcal B_{1}),(\frak g,[-,-]_{\frak g}),(\frak g^*,[-,-]_{\frak g^*})\big)$ such that the natural nondegenerate symmetric bilinear form $\mathcal{B}_{1}$ is invariant on $(\frak g\oplus \frak g^{*},[-,-]_{d})$.
\end{defi}

Recall that a {\bf Lie bialgebra} \cite{Cha,Drinfeld1983} is a triple $(\frak g,[-,-]_{\frak g},\delta)$ such that $(\frak g,[-,-]_{\frak g})$ is a Lie algebra, $(\frak g,\delta)$ is a Lie coalgebra and the following equation holds:
\begin{equation*}
\delta([x,y]_{\frak g})=\big( \mathrm{ad}_{\frak g}(x)\otimes\mathrm{id}_{\frak g}+\mathrm{id}_{\frak g}\otimes\mathrm{ad}_{\frak g}(x) \big)\delta(y)-\big( \mathrm{ad}_{\frak g}(y)\otimes\mathrm{id}_{\frak g}+\mathrm{id}_{\frak g}\otimes\mathrm{ad}_{\frak g}(y) \big)\delta(x),\;\forall x,y\in \frak g.
\end{equation*}
It is known that a Lie bialgebra is equivalent to a standard Manin triple of Lie algebras with respect to  the nondegenerate symmetric invariant bilinear form.

\begin{pro}\label{2920}
Let $(\frak g,[-,-]_{\frak g},\delta)$ be a Lie bialgebra. Then $( \frak g,[-,-]_{\frak g},[-,-]_{\frak g},\delta,\delta)$ is a self-module Lie bialgebra.
Conversely, suppose that $(\frak g,[-,-]_{\frak g},\circ_{\frak g},\delta,\theta)$ is a self-module Lie bialgebra.
If the following equations hold:
\begin{equation*}
[x,y]_{\frak g} =x\circ_{\frak g}y,\;
\delta(x) =\theta(x),\;\forall x,y \in\frak g,
\end{equation*}
then $(\frak g,[-,-]_{\frak g},\delta)$ is a Lie bialgebra.
\end{pro}

Let $(\frak g,[-,-]_{\frak g})$ be a Lie algebra and $r=\sum\limits_{i}a_{i}\otimes b_{i}\in\frak g\otimes \frak g$.
Set ${\bf C}(r)\in \frak g\otimes \frak g\otimes \frak g$ by
\begin{eqnarray*}
{\bf C}(r)=\sum_{i,j} [a_{i}, a_{j}]_{\frak g}\otimes b_{i}\otimes b_{j}+
a_{i}\otimes [b_{i},a_{j}]_{\frak g}\otimes b_{j}
+a_{i}\otimes a_{j}\otimes [b_{i},b_{j}]_{\frak g}.
\end{eqnarray*}
We say $r$ is a solution of the {\bf classical Yang-Baxter equation} ({\bf CYBE} in short) in $(\frak g,[-,-]_{\frak g})$ if ${\bf C}(r) =0$.
Moreover, by \cite{Cha}, an antisymmetric solution $r$ of the CYBE in a Lie algebra $(\frak g,[-,-]_{\frak g})$ gives rise to a Lie bialgebra $(\frak g,[-,-]_{\frak g},\delta_{r})$ which is called {\bf triangular}, where $\delta_{r}$ is given by
\begin{eqnarray}\label{eq:Lie cob}
\delta_{r}(x)=\big(\mathrm{ad}_{\frak g}(x)\otimes\mathrm{id}_{\frak g}+\mathrm{id}_{\frak g}\otimes\mathrm{ad}_{\frak g}(x)\big)r,\;\forall x\in  \frak g .
\end{eqnarray}
We have the following relationships between triangular Lie bialgebras and triangular self-module Lie bialgebras.

\begin{pro}\label{2480}
Let $(\frak g,[-,-]_{\frak g})$ be a Lie algebra and $r\in \frak g\otimes\frak g$ be an antisymmetric solution of the CYBE in  $(\frak g,[-,-]_{\frak g})$, which gives rise to a triangular Lie bialgebra $(\frak g,[-,-]_{\frak g},\delta_{r})$ through \eqref{eq:Lie cob}.
Then $r$ is also an antisymmetric  solution of the SLYBE in the self-module Lie algebra $(\frak g,[-,-]_{\frak g},[-,-]_{\frak g})$, which gives rise to a triangular self-module Lie bialgebra $(\frak g,[-,-]_{\frak g},[-,-]_{\frak g}$,
$\Delta_{r},\theta_{r})$
with respect to $\mathrm{id}_{\frak g}$.
Moreover, the comultiplications $\Delta_{r},\theta_{r}$ and $\delta_{r}$ satisfy the following relations:
\begin{equation*}\label{3896}
\Delta_{r}(x)=\theta_{r}(x)=\delta_{r}(x),\;\forall x\in\frak g.
\end{equation*}
\end{pro}

\subsection{Pre-Lie bialgebras as a classical example of self-module Lie bialgebras}\

We have the following relationships between pre-Lie algebras and self-module Lie algebras.

\begin{pro}\label{2515}
Let $(\frak g, \circ_{\frak g})$ be a pre-Lie algebra and $(\frak g,[-,-]_{\frak g})$ be the sub-adjacent Lie algebra of  $(\frak g, \circ_{\frak g})$ with $[-,-]_{\frak g}$ given by
\begin{equation}\label{eq:sub-adj-}
[x,y]_{\frak g}=x\circ_{\frak g}y-y\circ_{\frak g}x,\;\forall x,y\in\frak g.
\end{equation}
Then $(\frak g,[-,-]_{\frak g},\circ_{\frak g})$ is a self-module Lie algebra.
Conversely, let $(\frak g,[-,-]_{\frak g},\circ_{\frak g})$ be a self-module Lie algebra.
If \eqref{eq:sub-adj-} holds, then $(\frak g,\circ_{\frak g})$ is a pre-Lie algebra which is compatible with $(\frak g,[-,-]_{\frak g})$.
\end{pro}

Recall that a {\bf symplectic form} on a Lie algebra $(\frak g,[-,-]_{\frak g})$ is a nondegenerate antisymmetric bilinear form $\mathcal{B}$ on $(\frak g,[-,-]_{\frak g})$ which satisfies
\begin{eqnarray}\label{eq:2coc}
\mathcal{B}([x,y]_{\frak g},z)+\mathcal{B}([y,z]_{\frak g},x)+\mathcal{B}([z,x]_{\frak g},y)=0,\;\forall x,y,z\in\frak g.
\end{eqnarray}
By \cite{Chu,Ku1}, we have the following result.

\begin{pro}\label{2590}
Let $(\frak g,[-,-]_{\frak g} )$ be a Lie algebra, $\mathcal{B}$ be a nondegenerate bilinear form on $\frak g$ and $(\frak g,[-,-]_{\frak g}$,
$\circ_{\frak g})$ be the self-module Lie algebra induced from  $(\frak g,[-,-]_{\frak g},\mathcal{B} )$.
If $\mathcal{B}$ is a symplectic form on $(\frak g,[-,-]_{\frak g} )$, then $(\frak g,\circ_{\frak g})$ is a pre-Lie algebra which is compatible with $(\frak g,[-,-]_{\frak g})$.
\end{pro}

\begin{defi}\label{2388-}
A {\bf standard Manin triple of Lie algebras with respect to  the symplectic form} \cite{Bai2008} (or a {\bf parakähler Lie algebra}) is a standard generic Manin triple of Lie algebras  $\big( (\frak g\oplus\frak g^*,[-,-]_d,\mathcal B_{-1}),(\frak g,[-,-]_{\frak g}), (\frak g^{*},[-,-]_{\frak g^{*}})  \big)$  such that the natural nondegenerate antisymmetric bilinear form $\mathcal{B}_{-1}$ given by
    \begin{equation*}
        \mathcal{B}_{-1}(x+\xi ,y+\eta )=-\langle x,\eta \rangle+\langle \xi ,y\rangle,\;\forall x,y\in \frak g, \xi ,\eta \in \frak g^{*}
    \end{equation*}
    is a symplectic form on $(\frak g\oplus \frak g^{*},[-,-]_{d})$.
\end{defi}

Recall that a {\bf pre-Lie coalgebra} \cite{Bai2008}
is a pair $(\frak g,\alpha)$, where $\frak g$ is a vector space and
$\alpha:\frak g\rightarrow\frak g\otimes\frak g$
is a comultiplication whose linear dual defines a
pre-Lie algebra structure on $\frak g^*$.
A {\bf pre-Lie bialgebra} \cite{Bai2008}
is a triple $(\frak g,\circ_{\frak g},\alpha)$ consisting of
a pre-Lie algebra $(\frak g,\circ_{\frak g})$ and
a pre-Lie coalgebra $(\frak g,\alpha)$ satisfying
\begin{align*}
(\alpha-\tau\alpha)(x\circ_{\frak g}y)
&=(\mathrm{id}_{\frak g}^{\otimes 2}-\tau)
\Big(
\big(L_{\circ_{\frak g}}(x)\otimes
\mathrm{id}_{\frak g}\big)(\alpha-\tau\alpha)(y)
+\big(\mathrm{id}_{\frak g}\otimes
R_{\circ_{\frak g}}(y)\big)\alpha(x)
\Big),\\
\alpha(x\circ_{\frak g}y-y\circ_{\frak g}x)
&=\big(
L_{\circ_{\frak g}}(x)\otimes\mathrm{id}_{\frak g}
+\mathrm{id}_{\frak g}\otimes
(L_{\circ_{\frak g}}-R_{\circ_{\frak g}})(x)
\big)\alpha(y)\nonumber\\
&\quad-\big(
L_{\circ_{\frak g}}(y)\otimes\mathrm{id}_{\frak g}
+\mathrm{id}_{\frak g}\otimes
(L_{\circ_{\frak g}}-R_{\circ_{\frak g}})(y)
\big)\alpha(x),
\quad \forall x,y\in\frak g.
\end{align*}
It is known that a pre-Lie bialgebra is equivalent to a
standard Manin triple of Lie algebras with respect to  the symplectic form.

\begin{pro}\label{2646}
Let $(\frak g,\circ_{\frak g},\alpha)$ be a pre-Lie
bialgebra. Then
$(\frak g,[-,-]_{\frak g},\circ_{\frak g},\Delta,\alpha)$
is a self-module Lie bialgebra, where
$[-,-]_{\frak g}$ is given by \eqref{eq:sub-adj-} and
\begin{equation*}\label{3132}
\Delta(x)=\alpha(x)-\tau\alpha(x),\quad
\forall x\in\frak g.
\end{equation*}
Conversely, let
$(\frak g,[-,-]_{\frak g},\circ_{\frak g},\Delta,\theta)$
be a self-module Lie bialgebra.
If \eqref{eq:sub-adj-} and the following equation hold:
\begin{equation*}\label{3129}
\Delta(x)=\theta(x)-\tau\theta(x),\quad
\forall x\in\frak g,
\end{equation*}
then $(\frak g,\circ_{\frak g},\theta)$ is a pre-Lie
bialgebra.
\end{pro}

Let $(\frak g,\circ_{\frak g})$ be a pre-Lie algebra
and $r=\sum\limits_{i}a_{i}\otimes b_{i}\in \frak g\otimes \frak g$.
Set ${\bf S}(r)\in \frak g\otimes \frak g\otimes \frak g$ by
\begin{eqnarray*}
{\bf S}(r)=\sum_{i,j}-a_{i}\circ_{\frak g} a_{j}\otimes b_{i}\otimes b_{j}+a_{i}\otimes b_{i}\circ_{\frak g} a_{j}\otimes b_{j}+a_{i}\otimes a_{j}\otimes (b_{i}\circ_{\frak g}b_{j}-b_{j}\circ_{\frak g}b_{i}).
\end{eqnarray*}
We say $r$ is a solution of the {\bf S-equation} \cite{Bai2008} in $(\frak g,\circ_{\frak g})$ if ${\bf S}(r)=0$.
Moreover, by \cite{Bai2008},
a symmetric solution $r$ of the S-equation in a pre-Lie algebra $(\frak g,\circ_{\frak g})$ gives rise to a pre-Lie bialgebra $(\frak g,\circ_{\frak g},\alpha_{r})$ which is called {\bf triangular}, where $\alpha_{r}$ is given by
\begin{eqnarray}\label{eq:pre-Lie cob}
\alpha_{r}(x)=\big({L}_{\circ_{\frak g}}(x)\otimes\mathrm{id}_{\frak g}+\mathrm{id}_{\frak g}\otimes
(L_{\circ_{\frak g}}-R_{\circ_{\frak g}})(x)\big)r,\;\forall x\in \frak g.
\end{eqnarray}
We have the following relationships between triangular  pre-Lie bialgebras and triangular self-module Lie bialgebras.

\begin{pro}\label{2675}
Suppose that $(\frak g,\circ_{\frak g})$ is a pre-Lie algebra and $r\in\frak g\otimes\frak g$ is a symmetric solution of the $S$-equation in $(\frak g,\circ_{\frak g})$ which gives rise to a triangular pre-Lie bialgebra $(\frak g,\circ_{\frak g},\alpha_r)$ through \eqref{eq:pre-Lie cob}.
Then $r$ is also a symmetric solution of the SLYBE in the self-module Lie algebra $(\frak g,[-,-]_{\frak g},\circ_{\frak g})$ with $[-,-]_{\frak g}$ given by \eqref{eq:sub-adj-}, which gives rise to a triangular self-module Lie bialgebra $(\frak g,[-,-]_{\frak g},\circ_{\frak g},\Delta_r,\theta_r)$ with respect to  $-\mathrm{id}_{\frak g}$.
Moreover, the comultiplications $\Delta_{r},\theta_{r}$ and $\alpha_{r}$ satisfy the following relations:
\begin{align*}
\Delta_r(x) =\alpha_r(x)-\tau\alpha_r(x),\;
\theta_r(x) =\alpha_r(x),\;\forall x\in\frak g.
\end{align*}
\end{pro}

\subsection{Anti-pre-Lie bialgebras as a classical example of self-module Lie bialgebras}\

Let $\frak g$ be a vector space and $\star_{\frak g}:\frak g\otimes\frak g\rightarrow\frak g$ be a multiplication on $\frak g$.
Suppose that $(\frak g,\star_{\frak g})$ is a Lie-admissible algebra, that is,
$(\frak g,[-,-]_{\frak g})$ with $[-,-]_{\frak g}$ given by
\begin{equation}\label{eq:sub-adj2}
[x,y]_{\frak g}=x\star_{\frak g}y-y\star_{\frak g}x,\;\forall x,y\in\frak g
\end{equation}
is a Lie algebra. If the following equation holds:
\begin{equation*}
x\star_{\frak g}(y\star_{\frak g} z)-y\star_{\frak g}(x\star_{\frak g} z)
=(y\star_{\frak g}x-x\star_{\frak g}y)\star_{\frak g} z,\;\forall x,y,z\in\frak g,\label{eq:anti-pre-Lie1}
\end{equation*}
then we say $(\frak g,\star_{\frak g})$ is an {\bf anti-pre-Lie algebra} \cite{LB2022} which is {\bf compatible with $(\frak g,[-,-]_{\frak g})$}, and $(\frak g,[-,-]_{\frak g})$ is the {\bf sub-adjacent Lie algebra} of $(\frak g,\star_{\frak g})$.
Moreover, $(-L_{\star_{\frak g}},\frak g)$ is a representation of $(\frak g,[-,-]_{\frak g})$.

We have the following relations between anti-pre-Lie algebras and self-module Lie algebras.

\begin{pro}\label{2780}
Let $(\frak g,\star_{\frak g})$ be an anti-pre-Lie algebra.
Then there is a self-module Lie algebra $(\frak g,[-,-]_{\frak g},\circ_{\frak g})$, where
$[-,-]_{\frak g}$ is given by \eqref{eq:sub-adj2} and
$\circ_{\frak g}$ is given by
\begin{equation}\label{563}
x\circ_{\frak g}y=-x\star_{\frak g}y,\;\forall x,y\in\frak g.
\end{equation}
Conversely, let $(\frak g,[-,-]_{\frak g},\circ_{\frak g})$ be a self-module Lie algebra.
If  the following equation holds:
\begin{equation}\label{2907}
[x,y]_{\frak g}=y\circ_{\frak g}x-x\circ_{\frak g}y,\;\forall x,y\in\frak g,
\end{equation}
then $(\frak g,\star_{\frak g} )$ with $\star_{\frak g}$ given through \eqref{563}
is an anti-pre-Lie algebra.
\end{pro}

Recall  \cite{Dzh} that a {\bf  commutative $2$-cocycle}
on a Lie algebra $(\frak g,[-,-]_{\frak g})$ is a symmetric bilinear form $\mathcal{B}$ on $(\frak g,[-,-]_{\frak g})$ that satisfies \eqref{eq:2coc}.
By \cite{LB2022}, we have the following result.

\begin{pro}\label{2846}
Let $(\frak g,[-,-]_{\frak g} )$ be a Lie algebra, $\mathcal{B}$ be a nondegenerate bilinear form on $\frak g$ and $(\frak g,[-,-]_{\frak g}$,
$\circ_{\frak g})$ be the self-module Lie algebra induced from  $(\frak g,[-,-]_{\frak g},\mathcal{B} )$.
If $\mathcal{B}$ is a commutative $2$-cocycle on $(\frak g,[-,-]_{\frak g} )$, then $(\frak g,\star_{\frak g})$
with $\star_{\frak g}$ given through \eqref{563}
is an anti-pre-Lie algebra which is compatible with $(\frak g,[-,-]_{\frak g})$.
\end{pro}

\begin{defi}\label{2388--}
    A {\bf standard Manin triple of Lie algebras with respect to  the commutative $2$-cocycle} \cite{LB2023} is a  standard generic Manin triple of Lie algebras $\big((\frak g\oplus\frak g^*,[-,-]_d,\mathcal B_{1}),(\frak g,[-,-]_{\frak g}),(\frak g^*$,
    $[-,-]_{\frak g^*})\big)$ such that the natural nondegenerate symmetric bilinear form $\mathcal{B}_{1}$
    is a commutative $2$-cocycle on $(\frak g\oplus \frak g^{*},[-,-]_{d})$.
\end{defi}

Recall that an {\bf anti-pre-Lie coalgebra} \cite{LB2023}
is a pair $(\frak g,\beta)$, where $\frak g$ is a vector space and
$\beta:\frak g\rightarrow\frak g\otimes\frak g$
is a comultiplication whose linear dual defines an
anti-pre-Lie algebra structure on $\frak g^*$.
An {\bf anti-pre-Lie bialgebra} \cite{LB2023}
is a triple $(\frak g,\star_{\frak g},\beta)$ consisting of
an anti-pre-Lie algebra $(\frak g,\star_{\frak g})$ and
an anti-pre-Lie coalgebra $(\frak g,\beta)$ satisfying
\begin{align*}
(\beta-\tau\beta)(x\star_{\frak g}y)
&=(\mathrm{id}_{\frak g}^{\otimes 2}-\tau)
\Big(
\big(L_{\star_{\frak g}}(x)\otimes
\mathrm{id}_{\frak g}\big)(\beta-\tau\beta)(y)
-\big(\mathrm{id}_{\frak g}\otimes
R_{\star_{\frak g}}(y)\big)\beta(x)
\Big),\\
\beta(x\star_{\frak g}y-y\star_{\frak g}x)
&=\big(
\mathrm{id}_{\frak g}\otimes
(L_{\star_{\frak g}}-R_{\star_{\frak g}})(x)
-L_{\star_{\frak g}}(x)\otimes
\mathrm{id}_{\frak g}
\big)\beta(y)\\
&\quad-\big(
\mathrm{id}_{\frak g}\otimes
(L_{\star_{\frak g}}-R_{\star_{\frak g}})(y)
-L_{\star_{\frak g}}(y)\otimes
\mathrm{id}_{\frak g}
\big)\beta(x),
\quad \forall x,y\in\frak g.
\end{align*}
It is known that an anti-pre-Lie bialgebra is equivalent to a
standard Manin triple of Lie algebras with respect to  the commutative $2$-cocycle.

\begin{pro}\label{2932}
Let $(\frak g,\star_{\frak g},\beta)$ be an anti-pre-Lie
bialgebra. Then
$(\frak g,[-,-]_{\frak g},\circ_{\frak g},\Delta,\theta)$
is a self-module Lie bialgebra, where
$[-,-]_{\frak g}$ and $\circ_{\frak g}$ are given by
\eqref{eq:sub-adj2} and \eqref{563}, respectively, and
\begin{align}
\Delta(x)&=\beta(x)-\tau\beta(x), \nonumber\\
\theta(x)&=-\beta(x),\;
\forall x\in\frak g.\label{3171}
\end{align}
Conversely, let
$(\frak g,[-,-]_{\frak g},\circ_{\frak g},\Delta,\theta)$
be a self-module Lie bialgebra.
If \eqref{2907} holds and
\begin{equation*}\label{3157}
\Delta(x)=\tau\theta(x)-\theta(x),\;
\forall x\in\frak g,
\end{equation*}
then $(\frak g,\star_{\frak g},\beta)$ is an anti-pre-Lie
bialgebra, where $\star_{\frak g}$ and $\beta$ are respectively given through
\eqref{563} and \eqref{3171}.
\end{pro}

Let $(\frak g,\star_{\frak g})$ be an anti-pre-Lie algebra
and $r=\sum\limits_{i}a_{i}\otimes b_{i}\in \frak g\otimes \frak g$.
Set ${\bf T}(r)\in \frak g\otimes \frak g\otimes \frak g$ by
\begin{eqnarray*}
{\bf T}(r)= \sum_{i,j} a_{i}\star_{\frak g} a_{j}\otimes b_{i}\otimes b_{j}
+a_{i}\otimes b_{i}\star_{\frak g} a_{j}\otimes b_{j}
-a_{i}\otimes a_{j}\otimes (b_{i}\star_{\frak g}b_{j}-b_{j}\star_{\frak g}b_{i}).
\end{eqnarray*}
We say $r$ is a solution of the {\bf anti-pre-Lie Yang-Baxter equation \cite{LB2023}} ({\bf APL-YBE} in short) in $(\frak g,\star_{\frak g})$  if ${\bf T}(r)=0$.
By \cite{LB2023},
an antisymmetric solution $r$ of the APL-YBE in an anti-pre-Lie algebra $(\frak g,\star_{\frak g})$ gives rise to an anti-pre-Lie bialgebra $(\frak g,\star_{\frak g},\beta_{r})$ which is called {\bf triangular}, where $\beta_{r}$ is given by
\begin{eqnarray}\label{eq:anti-pre-Lie cob}
\beta_{r}(x)=({L}_{\star_{\frak g}}(x)\otimes\mathrm{id}_{\frak g}-\mathrm{id}_{\frak g}\otimes(L_{\star_{\frak g}}-R_{\star_{\frak g}})(x))r,\;\forall x\in \frak g.
\end{eqnarray}
We have the following relationships between triangular anti-pre-Lie bialgebras and triangular self-module Lie bialgebras.

\begin{pro}\label{2972}
Suppose that $(\frak g,\star_{\frak g})$ is an anti-pre-Lie algebra and $r$ is an antisymmetric solution of the APL-YBE in $(\frak g,\star_{\frak g})$ which gives rise to a triangular anti-pre-Lie bialgebra $(\frak g,\star_{\frak g},\beta_{r})$ through \eqref{eq:anti-pre-Lie cob}.
Then  $r$ is also an antisymmetric solution of the SLYBE in the self-module Lie algebra $(\frak g,[-,-]_{\frak g},\circ_{\frak g})$ given in Proposition \ref{2780},
which gives rise to a triangular self-module Lie bialgebra $(\frak g,[-,-]_{\frak g},\circ_{\frak g},\Delta_{r},\theta_{r})$
with respect to $\mathrm{id}_{\frak g}$.
Moreover, the comultiplications $\Delta_{r},\theta_{r}$ and $\beta_{r}$ satisfy the following relations:
\begin{align*}
\Delta_{r}(x) =\beta_{r}(x)-\tau\beta_{r}(x),\;
\theta_{r}(x) =-\beta_{r}(x),\;\forall x\in\frak g.
\end{align*}
\end{pro}
The three classical examples discussed in this section are
summarized in Table~\ref{tab:classical-specializations}.

\begin{table}[htbp]
    \centering
    \caption{Classical examples of self-module Lie bialgebras.}
    \label{tab:classical-specializations}
    \begingroup
    \footnotesize
    \renewcommand{\arraystretch}{1.7}
    \setlength{\tabcolsep}{3pt}
    \begin{tabular}{|p{0.18\textwidth}|p{0.23\textwidth}|p{0.24\textwidth}|p{0.26\textwidth}|}
        \hline
        & Lie & Pre-Lie & Anti-pre-Lie \\
        \hline
        Algebra relation
        & $x\circ_{\frak g}y=[x,y]_{\frak g}$
        & $[x,y]_{\frak g}=x\circ_{\frak g}y-y\circ_{\frak g}x$
        & $\begin{aligned}[t]
            &x\star_{\frak g}y=-x\circ_{\frak g}y,\\
            &[x,y]_{\frak g}=y\circ_{\frak g}x-x\circ_{\frak g}y
        \end{aligned}$ \\
        \hline
        Pairing on $\frak g\oplus\frak g^*$
        & $\mathcal B_{1}$:\newline
        symmetric invariant
        & $\mathcal B_{-1}$:\newline
        symplectic
        & $\mathcal B_{1}$:\newline
        commutative $2$-cocycle \\
        \hline
        Coalgebra relation
        & $\Delta=\theta=\delta$
        & $\Delta=\alpha-\tau\alpha,\;\theta=\alpha$
        & $\Delta=\beta-\tau\beta,\;\theta=-\beta$ \\
        \hline
        Tensor equation
        & CYBE
        & $S$-equation
        & APL-YBE \\
        \hline
        Tensor symmetry
        & Antisymmetric
        & Symmetric
        & Antisymmetric \\
        \hline
        Triangular relation
        & $\Delta_r=\theta_r=\delta_r$
        & $\Delta_r=\alpha_r-\tau\alpha_r,\;\theta_r=\alpha_r$
        & $\Delta_r=\beta_r-\tau\beta_r,\;\theta_r=-\beta_r$ \\
        \hline
    \end{tabular}
    \endgroup
\end{table}

\FloatBarrier
\begin{rmk}\label{rmk:c-ex}
The first self-module Lie bialgebra in Example \ref{ex:nonclassical} does not belong to any of the three classical specializations in Table~\ref{tab:classical-specializations}, since
\begin{eqnarray*}
    \Delta_s(e_2)\neq\theta_s(e_2),\;
    \Delta_s(e_2)\neq0=\theta_s(e_2)-\tau\theta_s(e_2).
\end{eqnarray*}
The second self-module Lie bialgebra in the same example recovers a triangular Lie bialgebra, since $\circ_2=[-,-]_{\frak h}$ and $\Delta_r=\theta_r$.
\end{rmk}

\noindent{\bf Author contributions. } The authors made equal contributions to the paper.

\noindent{\bf Acknowledgements.}
This work is supported by NSFC (12271265, 12401031, W2412041, 12671037),
the China Postdoctoral Science Foundation (2024T005TJ,
2024M761507), the Postdoctoral Fellowship Program of CPSF (GZC20240755), Central Government-Guided Fund Project for Local Science and Technology Development in Tianjin (25ZYZYJC00010),  the Fundamental Research Funds for the Central Universities and Nankai Zhide Foundation.

\noindent{\bf Declaration of interests.}
The authors have no conflicts of interest to disclose.

\noindent{\bf Data availability.}
Data sharing is not applicable as no data were created or analyzed.

\end{document}